\documentclass[11pt]{amsart}

\newtheorem{theorem}{Theorem}[section]
\newtheorem{lemma}[theorem]{Lemma}

\theoremstyle{definition}
\newtheorem{definition}[theorem]{Definition}

\theoremstyle{remark}
\newtheorem{remark}[theorem]{Remark}

\numberwithin{equation}{section}

\usepackage{a4,a4wide}
\usepackage{leqno}
\usepackage{graphics}
\usepackage{amssymb}
\usepackage{amsmath}
\usepackage{amsthm}
\usepackage[latin1]{inputenc}
\usepackage{euscript}
\usepackage{palatino}
\usepackage{color}
\newenvironment{theo*}[2]{\smallskip\newline \noindent\textbf{Theorem}\;(#1)\;\newline \textit{#2}\medskip\\}

\newenvironment{th*}[1]{\smallskip\newline \noindent \textbf{Theorem}\;--\;\textit{#1}\smallskip}

\newcommand{\beq}{\begin{equation}}
\newcommand{\eeq}[1]{\label{#1}\end{equation}}

\def\O{{\Omega}}
\def\o{{\omega}}

\def\eps{{\epsilon}}
\def\G{{\Gamma}}

\def\e{{\mathcal{E}}}

\def\m{{\mathcal{M}}}

\def\r{{\mathcal{R}}}

\def\t{{\mathcal{T}}}

\def\L{{\mathcal{L}}}

\def\S{\Sigma}
\def\M{{\mathcal{M}}}

\def\R{{\mathbb{R}}}

\def\N{{\mathbb{N}}}

\newcommand{\under}[1]{\underline{#1}}
\newcommand{\ope}[1]{\e[{#1}]}

\newcommand{\lb}[1]{\L_{_{#1}}}
\newcommand{\oplb}[2]{\L_{_{#2}}[{#1}]}
\newcommand{\opm}[1]{\M[{#1}]}
\newcommand{\mb}[1]{\m_{_{#1}}}
\newcommand{\opmb}[2]{\M_{_{#2}}[{#1}]}

\newcommand{\opr}[1]{\r[{#1}]}

\newcommand{\dem}[1]{\vskip 0.2\baselineskip \noindent {\bf{#1}}\vskip 0.2\baselineskip }
\newcommand{\fdem}{\vskip 0.2 pt \qquad \qquad \qquad \qquad \qquad \qquad \qquad \qquad \qquad \qquad \qquad \qquad \qquad \qquad \qquad \qquad  $\square$  }

\newtheorem{prop}{\textbf{Proposition --}}[section]
\newtheorem{cla}{\textbf{Claim --}}[section]

\def\tilde{\widetilde}

\begin{document}

\title{On a simple criterion for the existence of a principal eigenfunction of some nonlocal operators }

\author{J{\'e}r{\^o}me Coville\\
~\\
\textit{\tiny  UR 546 Biostatistique et Processus Spatiaux\\
 INRA, Domaine St Paul Site Agroparc\\
  F-84000 Avignon\\
   France}}
\address{UR 546 Biostatistique et Processus Spatiaux\\
 INRA, Domaine St Paul Site Agroparc\\
  F-84000 Avignon\\
 France\\}
\email{jerome.coville@avignon.inra.fr}
\thanks{The author is supported by INRA Avignon and  is thankful to the  Max Planck Institute for Mathematics in the Science where part of  this work has been done. The author want also to warmly thank the anonymous referee for the numerous useful comments  he has made to improve this paper.  }

\subjclass[2000]{ }




\begin{abstract}
In this paper we are interested in the existence of a principal eigenfunction of a nonlocal operator which appears in the description of various phenomena ranging from population dynamics to micro-magnetism. More precisely, we study the following eigenvalue problem:
$$\int_{\O}J\left(\frac{x-y}{g(y)}\right)\frac{\phi(y)}{g^n(y)}\, dy +a(x)\phi =\rho \phi,$$
where $\O\subset\R^n$ is an open connected set, $J$  a nonnegative kernel and $g$ a positive function. 
First, we establish  a criterion for the existence of a principal eigenpair $(\lambda_p,\phi_p)$. We also  explore the relation between the  sign of the largest element of the spectrum with a strong maximum property satisfied by the operator.  As an application of these results we construct and  characterize  the solutions of some nonlinear nonlocal reaction diffusion equations. 
\end{abstract}

\maketitle
\section{\bf Introduction and Main results}
In the past few years much attention has been drawn to the study of nonlocal reaction diffusion equations, where the usual elliptic diffusion operator is replaced by a nonlocal operator of the form
\begin{equation}
\opm{u}:=\int_{\O}k(x,y)u(y)\,dy -b(x)u, \label{pev-eq-gene}
\end{equation}
  where $\O\subset \R^n$, $k \ge 0$ satisfies  $\int_{\R^n}
k(y,x)dy<\infty$ for all $x\in \R^n$ and $b(x)\in C(\O)$; see among other references \cite{AB,BC1,BFRW,CCR,Ch,CR,CCEM,CER,Co2,Co4,CD1,CD2,DGP,KM,M,SSN}.
Such type of diffusion process has been widely used to describe the dispersal of a  population through its environment in the
following sense. As stated in  \cite{F1,F2,HMMV} if $u(y,t)$ is
thought of as a density at  a location $y$ at a time $t$ and $k(x,y)$
as the probability distribution of jumping from  a location $y$ to a
location $x$, then the rate at which the individuals from all other
places are arriving to the location $x$ is
$$\int_{\O} k(x,y)u(y,t)\,dy.$$
On the other hand, the rate at which the individuals are leaving the location $x$  is $-b(x)u(x,t)$.
This formulation of the dispersal of individuals finds its justification in many ecological problems of seed dispersion; see for example \cite{CMS,Cl,DK,KM,M,SSN}.

In this paper, we study the properties of the principal eigenvalue  of the operator $\m$, when the  kernel $k(x,y)$ takes the form  
\begin{equation}
k(x,y) = J\left(\frac{x-y}{g(y)}\right)\frac{1}{g^n(y)},\label{pev-eq-ccem}
\end{equation}
where $J$ is a continuous probability density and the function $g$  is  bounded and positive.
That is to say we investigate the following eigenvalue problem:
\begin{equation}\label{pev.eq.ev}
\int_{\O}J\left(\frac{x-y}{g(y)}\right)\frac{u(y)}{g^n(y)}\,dy -b(x)u=-\lambda u \quad \text{in} \quad \O.
\end{equation}
Such type of diffusion kernel was recently introduced by Cortazar \textit{et al.} \cite{CCEM} in order to model a  non homogeneous dispersal process.  
Along this paper, with no further specifications, we will always make the following assumptions on   $\O$, $J$, $g$ and $b$ : 
\begin{align*}
&\O\subset \R^n \quad \text{is an open connected set}&\qquad(H1) \\ 
&J \in C_c(\R^n),\, J\ge 0,\, J(0)>0 &\qquad(H2)\\
&g\in L^{\infty}(\O),\, 0<\alpha\le g \le \beta,\, &\qquad(H3)\\
&b \in C(\bar\O)\cap L^{\infty}(\O) &\qquad(H4)
\end{align*}
where $C_c(\R^n)$ denotes the set of continuous functions with compact support.

The existence and a variational  characterization of the principal eigenvalue $\lambda_p$ of $\m$ is known from a long time, see for example Donsker and Varadhan \cite{DV}. However, as Donsker and Varadhan \cite{DV} have already noticed,  $\lambda_p$ is in general not an eigenvalue, that is to say there exists no positive function $\phi_p$ such that $(\lambda_p,\phi_p)$ is a solution of \eqref{pev.eq.ev}.  In this paper, we are  interested in finding some conditions on $\m$ ensuring the existence of a principal eigenpair $(\lambda_p,\phi_p)$  of \eqref{pev.eq.ev} such that  $\phi_p\in C(\O)$ and $\phi_p>0$. Such type of solution is commonly used  to analyse the long-time behaviour of some nonlocal evolution problems \cite{CCR,CCEM} and had proven to be a very efficient tool  in the analysis of nonlinear integrodifferential problems; see for example \cite{CDM1,GR}. 

 To our knowledge,  besides some particular situations 
the existence of an  eigenpair $(\lambda_p,\phi_p)$ for the equation \eqref{pev.eq.ev} is still an open question and many of the known results  concern these two cases:
\begin{enumerate}
\item $b(x)\equiv Constant$ 
\item The operator $\m$ satisfies a mass preserving property, i.e  $\forall u\in C(\O)$,
$$\int_{\O}\int_{\O}J\left(\frac{x-y}{g(y)}\right)\frac{u(y)}{g^n(y)}\,dydx -\int_{\O} b(x)u(x)\,dx =0.$$ 
\end{enumerate}
In both cases, the principal eigenvalue problem \eqref{pev.eq.ev} is  either reduced to  the analysis of the spectrum of the positive operator $\lb{\O}$ defined below: 
 $$\oplb{u}{\O}:=\int_{\O}J\left(\frac{x-y}{g(y)}\right)\frac{u(y)}{g^n(y)}\,dy$$ or  the principal eigenvalue is explicitly known, i.e. $\lambda_p=0$ and the principal eigenfunction  $\phi_p$ is also the positive solution of the following eigenvalue problem
 $$ \int_{\O}J\left(\frac{x-y}{g(y)}\right)\frac{\psi(y)}{g^n(y)}\,dy =\rho b(x)\psi.$$ 
Note, that even in this two simplified cases, showing the existence of an eigenfunction is still a difficult task when the domain $\O$ is unbounded.   

As observed in \cite{Co4},  the equation \eqref{pev-eq-gene} shares many properties with the usual  elliptic operators 
$$\e:=\sigma_{ij}(x)\partial_{ij} +\beta_i(x)\partial_i +c(x).$$ In particular, acting on smooth functions, we can rewrite $\m$
$$\opm{u}=\ope{u}+ \opr{u}$$ with  $\r$  an operator involving derivatives of higher order that in $\e$. 

Indeed, we have 
$$
\opm{u}=\int_{\O}k(x,y)[u(y)-u(x)]\, dy -c(x)u,
$$
with $c(x):=b(x)-\int_{\O}k(x,y)dy$. 
Using the change of variables  $z=x-y$ and performing a formal Taylor expansion of $u$ in the integral,  we can rewrite the nonlocal operator as follows   
$$\int_{x-\O}k(x,x-z)[u(x-z)-u(x)]\, dy=\sigma_{ij}(x)\partial_{ij}u +\beta_i(x) \partial_iu + \opr{u}$$
where we use the Einstein summation convention and $\sigma_{ij}(x)$, $\beta_i(x)$, and $\r$ are defined by the following expressions  
\begin{align*}
&\sigma_{ij}(x)=\frac{1}{2}\int_{x-\O}k(x,x-z)z_iz_j\,dz\\
&\beta_{i}(x)=\int_{x-\O}k(x,x-z)z_i\,dz \\
&\opr{u}:=\int_{0}^{1}\int_{0}^{1}\int_{0}^{1}\int_{x-\O}k(x,x-z)z_iz_jt^2s\partial_{ijk}u(x+ts\tau z) \,dtdsd\tau dz.
\end{align*}

For a second order elliptic operator $\e$,  the existence of a principal eigenpair $(\lambda_p,\phi_p)$ is well known  and  various  variational formulas characterising the principal eigenvalue exist, see for example \cite{BNV,DV,E, NP, PW,P}. In particular,  Berestycki, Nirenberg and Varadhan \cite{BNV} give a very simple and general definition of the principal eigenvalue of $\e$  that we recall below. 
Namely,  they define the principal eigenvalue of the elliptic operator $\e$  by the following quantity:
\begin{equation}
\label{pev.ellip.def}
\lambda_1:=\sup\left\{\lambda\in \R\,| \exists \phi \in\, C(\O), \phi>0,\; \text{ such that }\; \ope{\phi} +\lambda\phi \le 0\right\}.
\end{equation}

In this paper, we  adopt the definition of Berestycki, Nirenberg and Varadhan for the  definition of  the principal eigenvalue of the operator $\m$. The principal eigenvalue of the operator $\m$ is then given by the following quantity:  
 \begin{equation*}
\lambda_p(\m):=\sup\left\{\lambda\in \R\,| \exists \phi \in\, C(\O), \phi>0,\; \text{ such that }\; \opm{\phi} +\lambda\phi \le 0\right\}.
\end{equation*}

To make more explicit the dependence of the different parameters and to simplify the presentation of the results, we shall adopt the following notations:
\begin{itemize}
\item $a(x):=-b(x)$
\item $\sigma:=sup_{\O}a(x)$
\item $d\mu$ is the measure defined by $d\mu:=\frac{dx}{g^{n}(x)}$
\item $\oplb{u}{\O}:=\int_{\O}J\left(\frac{x-y}{g(y)}\right)\frac{u(y)}{g^{n}(y)}\,dy,=\int_{\O}J\left(\frac{x-y}{g(y)}\right)u(y)\,d\mu$
\item $\m :=\mb{\O}:=\lb{\O} + a(x)Id$
\end{itemize}
With this new notation the principal eigenvalue of $\mb{\O}$ can be rewritten as follows
\begin{equation}
\label{pev.eq.def} 
\lambda_p(\mb{\O}):=\sup\left\{\lambda\in \R\,| \exists \phi \in\, C(\O), \phi > 0,\; \text{ such that }\; \oplb{\phi}{\O} +(a(x)+\lambda)\phi \le 0\right\}.
\end{equation}

 Under the assumptions $(H1-H4)$, the principal eigenvalue $\lambda_p(\mb{\O})$ is well defined, see the appendix for the details. 

Obviously,  $\lambda_p$ is  monotone  with respect to the domain, the zero order term $a(x)$ and $J$. Moreover, $\lambda_p$ is a concave function of its argument and is Lipschitz continuous with respect to  $a(x)$. More precisely, we have 

\begin{prop}~\label{pev.prop1}
\begin{itemize}
\item[(i)] Assume $\O_1\subset\O_2$, then $$
\lambda_p(\lb{\O_1}+a(x))\ge \lambda_p(\lb{\O_2}+a(x)).
$$
\item[(ii)]Fix $\O$  and assume that $a_1(x)\ge a_2(x)$, then 
$$
\lambda_p(\lb{\O}+a_2(x))\ge\lambda_p(\lb{\O}+a_1(x)).
$$
Moreover, if $a_1(x)\ge a_2(x)+\delta$ for some $\delta>0$ then 
$$
\lambda_p(\lb{\O}+a_2(x))> \lambda_p(\lb{\O}+a_1(x)).
$$
\item[(iii)] $\lambda_p(\lb{\O}+a(x))$ is Lipschitz continuous in $a(x)$. More precisely,
$$|\lambda_p(\lb{\O}+a(x))- \lambda_p(\lb{\O}+b(x))|\le \|a(x)-b(x)\|_{\infty}$$
\item[(iv)] Let $J_1\le J_2$ be two positive continuous integrable functions and let us denote respectively $\lb{1,\O}$ and $\lb{2,\O}$ the corresponding operators. Then we have   
$$
\lambda_p(\lb{1,\O}+a(x))>\lambda_p(\lb{2,\O}+a(x)).
$$
\end{itemize}
\end{prop}


Let us  state our  first result concerning a sufficient condition for the existence of a principal eigenpair $(\lambda_p,\phi_p)$ for the operator  $\m$.
\begin{theorem}[Sufficient condition] \label{pev.th1}
Assume that $\O$, $J$, $g$ and $a$ satisfy  $(H1-H4)$. Let us denote $\sigma:=\sup_{\bar \O} a(x)$ and assume further that the function   $a(x)$ satisfies $\frac{1}{\sigma -a(x)}\not \in L^1_{d\mu}(\O_0)$ for some bounded domain $\O_0\subset\bar \O$. Then there exists a principal eigenpair $(\lambda_p,\phi_p)$ solution of \eqref{pev.eq.ev}. Moreover, $\phi_p\in C(\O)$,  $\phi_p>0$ and we have the following estimate
$$ -\sigma' <\lambda_p<-\sigma, $$ where $\sigma':=\sup_{x\in\O } \left[a(x)+ \int_{\O}J\left(\frac{y-x}{g(x)}\right)\frac{dy}{g^{n}(x)}\right]$.
\end{theorem} 

Note  that the Theorem holds true whenever  $\O$ is bounded or not.

The condition  $\frac{1}{\sigma -a(x)}\not \in L^1_{d\mu}(\O_0)$ is sharp in the sense that  if  $\frac{1}{\sigma -a(x)} \in L^1_{d\mu,loc}(\O)$ then we can construct  an operator $\mb{\O}$ such that the equation \eqref{pev.eq.ev} does not have a principal eigenpair. This is discussed in section \ref{pev.s.ce}, where  such an operator is constructed. 
We want also to stress that the boundedness of the open set $\O$ does not ensure the existence of an eigenfunction, see the counterexample in section \ref{pev.s.ce}.

In contrast with the elliptic case,  the sufficient condition has nothing to do with the regularity of the functions $a(x)$, $J$ or $g$.  This means that in general improving the regularity of the coefficients does not ensure at all the existence of an eigenpair. However, in low dimension of space $n=1,2$ the condition $\frac{1}{\sigma-a(x)} \not \in L^1_{d\mu}(\O_0)$  can be related to a regularity condition on the coefficient $a(x)$. Indeed, in one dimension if $a$ is Lipschitz continuous and achieves a maximum in $\O$ then the condition $\frac{1}{\sigma-a(x)} \not \in L^1_{d\mu}(\O_0)$ is automatically satisfied. Similarly, when $n=2$ the non-integrability condition is always satisfied when $a(x) \in C^{1,1}(\O)$ and achieves a maximum in $\O$.  More precisely, we have the following
\begin{theorem}
Assume that $\O$, $J$, $g$ and $a$ satisfy  $(H1-H4)$,  that $a$ achieves a global maximum at some point $x_0\in \O$. Then there exists a principal eigenpair $(\lambda_p,\phi_p)$ solution of \eqref{pev.eq.ev}  in the following situations   
\begin{itemize}
\item  $n=1, a(x) \in C^{0,1}(\O)$
\item  $n=2, a(x) \in C^{1,1}(\O) $
\item  $n\ge 3, a(x)\in C^{n-1,1}(\O), \forall k<n,  \partial^{k}a(x_0)=0$. 
\end{itemize}
\end{theorem}

One of the most interesting properties of the principal eigenvalue for an elliptic operator $\e$ is its relation with the existence of a maximum principle for $\e$. Indeed,   Berestycki \textit{ et al.} \cite{BNV}  have shown that there exists a strong  relation between the sign of this principal eigenvalue and the existence of a maximum principle for the elliptic operator $\e$.
Namely, they have proved  
\begin{theorem}[BNV]
Let $\O$ be a bounded open set, then $\e$ satisfies a refined maximum principle if and only if $\lambda_1>0$.
\end{theorem}

It turns out that when the principal eigenpair exists for  $\m$, we can  also obtain a similar relation between  the sign of the principal eigenvalue  of $\m$ and  some  maximum principle property.
More precisely, let us first define the maximum principle property satisfied by $\m$:
\begin{definition}[Maximum principle]\label{pev.def.pm}
When $\O$ is bounded, we say that the maximum principle is satisfied by  an operator $\mb{\O}$  if  for all function $u\in C(\bar\O)$ satisfying
\begin{align*}
&\opmb{u}{\O}\le 0  \quad\text{ in }\quad \O\\
&u\ge 0  \quad\text{ in }\quad \partial\O\
\end{align*}
then $u\ge 0$ in $\O$.
\end{definition}
 
With this definition of maximum principle,    we   show    

 \begin{theorem}\label{pev.th2}
 Assume that  $\O$ is a bounded set and let $J$, $g$ and $a$ be as in Theorem \ref{pev.th1}. Then the maximum principle is satisfied  by $\mb{\O}$ if and only if  $\lambda_p(\mb{\O})\ge 0$.  
 \end{theorem}

Note that there is a slight difference between the criteria for elliptic operators and for nonlocal ones. To have a maximum principle  for nonlocal operator it is sufficient to have a non negative principal eigenvalue, which is untrue for a elliptic operator where a strict sign of $\lambda_p$ is required.

Our last result is  an application of the sufficient condition for the existence of a principal  eigenpair to obtain a simple criterion for the  existence/non-existence  of a positive solution of the following  semilinear problem:
\begin{equation}
\opmb{u}{\O} +f(x,u)=0 \quad \text{ in }\quad \O.\label{pev.eq.semilin}
\end{equation}
where $f$ is a KPP type non-linearity. Such type of equation naturally appears in some ecological problems when in addition to the dispersion of the individuals in the environment,  the birth and death of these individuals  are also modelled, see \cite{GR,HMMV,KM, M}.

On $f$ we assume that:
\begin{align}
\label{hyp f1} \left\{
\begin{aligned}
& \hbox{$f \in C(\R\times[0,\infty))$ and is differentiable with respect to $u$} \\
& \hbox{$f_u(\cdot,0)$ is Lipschitz}  \\
& \hbox{$f(\cdot,0)\equiv 0$ and $f(x,u)/u$ is decreasing with
respect to $u$}
\\
& \hbox{there exists $M>0$ such that $f(x,u)\le0$ for all $u \ge
M$ and all $x$.}
\end{aligned}
\right.
\end{align}
The simplest example of such a nonlinearity is
$$
f(x,u) = u (\mu(x) - u ),
$$
where $\mu(x)$ is a Lipschitz function.

Such type of problem have received  recently a lot of attention, see for example  \cite{BZ,HMMV,KM,M} and reference therein. In particular, for $\O$  bounded and for a symmetric kernel $J$ Hutson \textit{et al.}  \cite{HMMV} have shown that there exists  a unique non trivial stationary solution  \eqref{pev.eq.semilin} provided that some principal eigenvalue of the linearised operator around the solution $0$ is positive.   
This result can be extended to more general kernel $J$ using the definition of principal eigenvalue \eqref{pev.eq.def}. More precisely, we show that 

\begin{theorem}
\label{pev.th3} Assume $\O$, $J$, $g$ and $a$  satisfy (H1-H4), $\O$ is bounded, $a(x)\le 0$ and $f$
satisfies \eqref{hyp f1}. Then there exists a unique non trivial
 solution of \eqref{pev.eq.semilin} when
$$
\lambda_p(\mb{\O}+f_u(x,0))<0,
$$
where $\lambda_p$ is the principal eigenvalue of the linear
operator {$\mb{\O}+f_u(x,0)$}. 
Moreover,  if   $\lambda_p \ge 0$ then any  nonnegative uniformly bounded solution  of \eqref{pev.eq.semilin} is identically zero.
\end{theorem}

As a consequence, we can derive the asymptotic behaviour of the solution of  the evolution  problem associated to  \eqref{pev.eq.semilin}:
\begin{align}
&\frac{\partial u}{\partial t}=\opmb{u}{\O} +f(x,u) \quad \text{ in }\quad \R^+\times\O.\label{pev.eq.semilin-para}\\
&u(0,x)=u_0(x)\quad \text{ in }\quad \O
\end{align}
Namely,  the asymptotic behaviour of $u(t, x)$ as $t \to +\infty$ is described in the following theorem:
\begin{theorem}\label{pev.thab}
Let $\O, J, g, b$ and $f$  be as in Theorem \ref{pev.th3}. Let $u_0$ be an arbitrary
bounded and continuous function in $\O$ such that $u_0 \ge 0, u_0 \not\equiv 0$. Let $u(t, x)$
be the solution of \eqref{pev.eq.semilin-para} with initial datum $u(0, x) = u_0(x)$. Then, we have
\begin{enumerate}
\item  If 0 is an unstable solution of \eqref{pev.eq.semilin} (that is $\lambda_p < 0$), then $u(t, x) \to p(x)$ pointwise as $t\to \infty$ , where $p$ is the unique positive solution of \eqref{pev.eq.semilin} given
by Theorem \ref{pev.th3}.
\item If $0$ is a stable solution of \eqref{pev.eq.semilin} (that is $\lambda_p\ge 0$), then $u(t, x)  \to 0$ pointwise
in $\O$ as $t\to +\infty$.
\end{enumerate}
\end{theorem}

Note that this criterion involves only the sign of $\lambda_p$ and  does not require any conditions on the function $f_u(x,0)$  ensuring the existence of a principal  eigenfunction.  Therefore, even in a situation where no principal eigenfunction exists for the operator $\mb{\O}+f_{u}(x,0)$ we still have  information on the survival or the extinction of the considered species. 
 Observe also that the condition obtained on the principal eigenvalue of the linearised operator is sufficient and necessary for the existence of a non trivial solution.

Before going into the proofs of these results, let us make some comments.
We first point out that the proofs we have given apply to a more general situation. More precisely, the above results can be easily  extended  to the case of  a dispersal kernel $k(x,y)$ which satisfies the following conditions:
\begin{align*}
&k(x,y) \in C_c(\O\times\O),\, k\ge 0,\, \int_{\O}k(x,y)\,dy<+\infty\quad \forall\, x\in \O &(\tilde H1)\\
&\exists\, c_0>0,\ \eps_0 >0\; \text{ such that } \min_{x\in\O}\left(\min_{y\in B(x,\eps_0)}k(x,y)\right)>c_0.&(\tilde H2)
\end{align*}
An example of such kernel is given by  
$$k(x,y)=J\left(\frac{x_1-y_1}{g_1(y)};\frac{x_2-y_2}{g_2(y)};\ldots; \frac{x_n-y_n}{g_n(y)}\right)\frac{1}{\prod_{i=1}^ng_i(y)},$$ 
with $0<\alpha_i\le g_i\le \beta_i$.

We want also to emphasize that the  condition that $J$ or $k$ has  a compact support is only needed to construct an eigenpair when $\O$ is unbounded. For a bounded domain, all the results  will also holds true if  $J$ is not assume compactly supported in $\O$.  

Note  that  the assumption $J(0)>0$ implies that the operator $\lb{\O}$ is not trivial on any open subset $\o\subset \O$, i.e. $\forall\, \o\subset\O,\;\forall\, u\in C(\O),\, \oplb{u}{\O}\neq 0 $ for $x\in \o$. This condition makes sure that the principal eigenfunction $\phi_p$ is positive in $\O$, which is a necessary condition for the existence of such principal eigenfunction.
Indeed, when  there exists an open subset $\o \subset \O$ such that $\lb{\O}$ is trivial, there is no guarantee that a principal eigenpair exists . For example,  this is the case for the operator $\mb{\O}$ where $\O:=(-1,1)$, $,J$ is such that $supp(J)\subset (\frac{1}{2}, 1)$ and  $3\le g\le 4$. In this situation, we easily see that for any $x \in (-\frac{1}{4},\frac{1}{4})$ and for any function $u \in C(\O)$,  we have $\oplb{u}{\O}(x)=0$. Therefore, the existence of an eigenfunction will strongly depend on the behaviour of the function $a(x)$ on this subset, i.e. $(\lambda_p+a(x))\phi\equiv 0 $ for $x\in (-\frac{1}{4},\frac{1}{4})$.  
If $(\lambda_p+a(x))\neq 0$  then $\phi\equiv 0$ in  $(-\frac{1}{4},\frac{1}{4})$.
 In this situation there is clearly  no existence of a positive principal eigenfunction.  
However, the condition $J(0)>0$ can still be relaxed and the above Theorems hold also true if we only assume that the kernel $J$ is such  that  there exists a positive integer  $p\in \N_{0}$   such that  the following kernel $J_p(x,y)$ satisfies $(\tilde H2)$ 
where $J_p(x,y)$ is defined by the recursion
\begin{align*}
&J_1(x,y):=J\left(\frac{x-y}{g(y)}\right)\frac{1}{g^{n}(y)}\\
&J_{p+1}(x,y):=\int_{\O} J_{p}(x,z)J_1(z,y)\,dz \quad \text{for } \quad p\ge 1.
\end{align*}

The above  condition is slightly more general that $J(0)>0$  and we see that $J(0)>0$ implies that $J_1$ satisfies  $(\tilde H2)$. 
In particular, as showed for example in \cite{AC}, for a convolution operator $K(x,y):=J(x-y)$, this new condition is optimal and can be related to a geometric condition on the convex hull of $\{y\in\R^n| J(y)>0\}$: 
\smallskip

\textit{There exists $p\in\N^*$, such that  $J_p$ satisfies $\tilde H2$   if and only if the convex hull of $\{y\in\R^n| J(y)>0\}$ contains $0$.}
 \medskip

We also want to stress that  we can easily extend the results of Theorems \ref{pev.th3} and \ref{pev.thab} to a periodic setting using the above generalisation on general nonnegative kernel.   Namely,  if we consider  the following problem   
\begin{equation}
\frac{\partial u}{\partial t}=\opmb{u}{\R^n} +f(x,u) \quad \text{ in }\quad \R^n\times\R^+,\label{pev.eq.semilinper}
\end{equation}
where $g$ and f(.,u) are assumed to be periodic functions  then the  existence of a unique non trivial
 periodic solution of \eqref{pev.eq.semilinper} is uniquely conditionned by the sign of the periodic principal eigenvalue 
$
\lambda_{p,per}(\mb{\R^n}+f_u(x,0)),
$
where $\lambda_{p,per}$ is defined as follows:
$$\lambda_{p,per}(\m):=\sup \{\lambda \in \R \,|\, \exists \psi>0, \psi\in C_{per}(\R^n) \text{ such that }\; \opmb{\psi}{\R^n}+\lambda\psi\le0\}.$$   
It is worth noticing that in this context, using  the periodicity, we have $$\lambda_{p,per}(\mb{\R^n}+f_u(x,0))=\lambda_{p}(\lb{Q}+f_u(x,0),Q),$$ 
 where $Q$ is the unit periodic cell and $\oplb{\psi}{Q}:=\int_{Q}k(x,y)u(y)dy$  with $k$ a positive kernel satisfying $\tilde H1$ and $\tilde H2$. Hence the analysis of the existence/ non existence of stationary solutions of \eqref{pev.eq.semilinper} will be handled through the the analysis of the existence/ non existence of stationary solutions of a  semilinear KPP  problem defined on a bounded domain.
\smallskip

Finally, along our analysis, provided a more restrictive assumption on the coefficient $a(x)$ is made,  we also observe that  Theorem \ref{pev.th1} holds as well  when we relax the assumption on the function $g$ and allow $g$ to touch $0$. More precisely, assuming that $g$ satisfies
$$
g\in L^{\infty}(\O),\, 0\le g \le \beta,\, \frac{1}{g^n} \in L^p_{loc}(\bar \O)\; \text{ with } \; p>1  \qquad(\tilde H3)
$$
then for a bounded domain $\O$, we have the following result:
\begin{theorem} \label{pev.th5}
Assume  that $\O$, $J$ and $a$ satisfy $(H1,H2,\tilde H3,H4)$, $\O$ bounded and $g$ satisfies $\tilde H3$. Let us denote $\sigma:=\sup_{\bar \O}a(x)$ and let $\G$  be the following set  
$$\G:=\{x\in \bar \O|\, a(x)=\sigma \}.$$
Assume further that  $\stackrel{\circ}{\G}\neq \emptyset$.   Then there exists a principal eigenpair $(\lambda_p,\phi_p)$ solution of \eqref{pev.eq.ev}. Moreover,  $\phi_p\in C(\O)$, $\phi_p>0$ and we have the following estimate
$$ -\sigma' <\lambda_p<-\sigma, $$ where $\sigma':=\sup_{x\in\O } \left[a(x)+ \int_{\O}J\left(\frac{y-x}{g(x)}\right)\frac{dy}{g^{n}(x)}\right]$.
\end{theorem}

 As a consequence the  criterion on the survival/extinction of a species obtained in Theorems \ref{pev.th3} and \ref{pev.thab} can be extended to such type of dispersal kernel.  
  More precisely, we have 
  \begin{theorem}
\label{pev.th6} Assume $\O$, $J$ and $g$  satisfy (H1,$\tilde H2,\tilde H3$), $\O$ is bounded  and $f$
satisfies \eqref{hyp f1}. Then there exists a unique non-trivial 
 solution of \eqref{pev.eq.semilin} if
$$
\lambda_p(\mb{\O}+f_u(x,0))<0,
$$
where $\lambda_p$ is the principal eigenvalue of the linear
operator {$\mb{\O}+f_u(x,0)$}. 
Moreover,  if   $\lambda_p \ge 0$ then any nonnegative uniformly  bounded solution is identically zero.
\end{theorem}
and 
\begin{theorem}\label{pev.thab2}
Let $\O, J, g, b$ and $f$  be as in Theorem \ref{pev.th6}. Let $u_0$ be an arbitrary
bounded and continuous function in $\O$ such that $u_0 \ge 0, u_0 \not\equiv 0$. Let $u(t, x)$
be the solution of \eqref{pev.eq.semilin-para} with initial datum $u(0, x) = u_0(x)$. Then, we have
\begin{enumerate}
\item  If 0 is an unstable solution of \eqref{pev.eq.semilin} (that is $\lambda_p < 0$), then $u(t, x) \to p(x)$ pointwise as $t\to \infty$ , where $p$ is the unique positive solution of \eqref{pev.eq.semilin} given
by Theorem \ref{pev.th6}.
\item If $0$ is a stable solution of \eqref{pev.eq.semilin} (that is $\lambda_p\ge 0$), then $u(t, x)  \to 0$ pointwise
in $\O$ as $t\to +\infty$.
\end{enumerate}
\end{theorem}

In this context, the existence of  a simple sufficient condition for the  existence of a principal eigenpair when $\O$ is an unbounded domain is more involved and  we have to make  a technical assumption on the set $\S:=\{x \in \bar\O,|\, g(x)=0\,\}$.    More precisely, we show 
 \begin{theorem} \label{pev.th7}
Assume  that $\O$, $J$ and $a$ satisfy $(H1,\tilde H2,H4)$ and $g$ satisfies $\tilde H3$. Let us denote $\sigma:=\sup_{\bar \O}a(x)$ and let $\G, \S$  be the following sets  
\begin{align*}
&\G:=\{x\in \bar \O|\, a(x)=\sigma \}\\
&\S:=\{x\in \bar \O|\, g(x)=0 \}.
\end{align*}
Assume further that $\O\cap\Sigma\subset\subset \O$   and   $\stackrel{\circ}{\G}\neq \emptyset$.   Then there exists a principal eigenpair $(\lambda_p,\phi_p)$ solution of \eqref{pev.eq.ev}. Moreover,  $\phi_p>0$ and we have the following estimate
$$ -\sigma' <\lambda_p<-\sigma, $$ where $\sigma':=\sup_{x\in\O } \left[a(x)+ \int_{\O}J\left(\frac{y-x}{g(x)}\right)\frac{dy}{g^{n}(x)}\right]$.
\end{theorem}

The paper is organized  as follows. In Section~\ref{pev.s.spectral} we review some spectral theory of positive operators and
we recall some Harnack's inequalities satisfied by a positive solution of  integral equation. Then, we prove the Theorems \ref{pev.th1} and  \ref{pev.th5} in Section \ref{pev.s.ev}.   The relation between the maximum principle  and the sign of the principal eigenvalue (Theorem  \ref{pev.th2})  and  a counter example to the existence of a principal eigenpair  are obtained respectively  in Section \ref{pev.s.mp} and in Section \ref{pev.s.ce}. The last two sections is devoted to the derivation of the survival/extinction criteria (Theorems \ref{pev.th3},\ref{pev.thab}, \ref{pev.th6}).

\section{Preliminaries \label{pev.s.spectral}}
In this section we first recall some results on the spectral theory of positive operators and  some Harnack's Inequalities satisfied by  a positive solution of 
\begin{equation}\label{pev.eq.har}
\oplb{u}{\O}-b(x)u=0,
\end{equation}
where $\lb{\O}$ is defined as above and   $b(x)$ is a positive continuous function in $ \O$.
Let us start with the spectral theory.
\subsection{Spectral Theory of positive operators}~\\
Let us recall some basic spectral results for positive
operators due to Edmunds, Potter and Stuart \cite{EPS} which are
extensions of the Krein-Rutman theorem for positive non-compact
operators.

A cone in a real Banach space $X$ is a non-empty closed set $K$
such that for all $x, y \in K$ and all $\alpha \ge 0$ one has $x +
\alpha y \in K$, and if $x\in K$, $-x \in K$ then $x=0$. A cone
$K$ is called reproducing if $X = K-K$. A cone $K$ induces a
partial ordering in $X$ by the relation $x\le y$ if and only if
$x-y \in K$. A linear map or operator $T:X\to X$ is called
positive if $T(K)\subseteq K$. The dual cone $K^*$ is the set of
functional $x^* \in X^*$ which are positive, that is, such that
$x^*(K) \subset [0,\infty)$.

If $T:X\to X$ is a bounded linear map on a complex Banach space X,
its essential spectrum (according to Browder \cite{browder})
consists of those $\lambda$ in the spectrum of $T$ such that at
least one of the following conditions holds : (1) the range of
$\lambda I - T$ is not closed, (2) $\lambda$ is a limit point of
the spectrum of $A$, (3) $\cup_{n=1}^\infty ker( (\lambda I -
T)^n)$ is infinite dimensional. The radius of the essential
spectrum of $T$, denoted by $r_e(T)$, is the largest value of
$|\lambda|$ with $\lambda$ in the essential spectrum of $T$. For
more properties of $r_e(T)$ see \cite{nussbaum}.

\begin{theorem}[Edmunds, Potter, Stuart]\label{pev.th.eps}~\\
Let K be a reproducing cone in a real Banach space X, and let
$T\in\L(X)$ be a positive operator such that $T^p(u)\ge cu$ for
some $u\in K$ with $\|u\|=1$, some positive integer $p$ and some
positive number $c$. Then if $c^{\frac{1}{p}}>r_e(T_c)$, $T$ has
an eigenvector $v\in K$ with associated  eigenvalue $\rho\ge
c^{\frac{1}{p}}$ and  $T^*$ has   eigenvector $v^*\in K^*$
corresponding to the eigenvalue  $\rho$. Moreover $\rho$ is
unique.
\end{theorem}
\medskip
A proof of this Theorem can be found in \cite{EPS}.


\subsection{Harnack's Inequality}~\\
Let us now present  some Harnack's  inequality satisfied by any positive continuous solution of the nonlocal equation \eqref{pev.eq.har} .


\begin{theorem}[Harnack Inequality]\label{pev.th.harnack}
Assume that $\O,J,g$ and $b>0$ satisfy $(H1,\tilde H2,H3,H4)$.
Let $\o\subset \subset \O$ be a compact set. Then there exists $C(J,\o,b,g)$ such that for all positive continuous bounded solutions $u$ of \eqref{pev.eq.har} we have 
$$ u(x)\le C u(y) \quad \text{ for all }\quad x,y \in \o.$$
\end{theorem}

When the assumption on $g$ is relaxed the above Harnack's estimate does not hold any more but an uniform estimate still holds.
Namely,  
\begin{theorem}[Local uniform  estimate]\label{pev.th.uniformesti}
 Assume that $\O,J,g$ and $b>0$ satisfy $(H1,\tilde H2,\tilde H3,H4)$. 
Assume that $\O\cap\Sigma\subset\subset \O$ and let $\o\subset\bar \O$ be a compact set. Let $\O(\o)$ denote the following set 
$$\O(\o):=\bigcup_{x\in \o}B(x,\beta).$$ 
Then there exists a positive constant $\eta^*$ such that, for any $0<\eta\le\eta^*$, there exist a compact set $\o'\subset \subset \O(\o)\cap \O$ and a constant $C(J,\o,\O,\o',b,g,\eta)$ such that the following assertions are verified
\begin{itemize} 
\item[(i)] $\{x\in \O(\o)\cap W_\eta|d(x,\partial(\O(\o)\cap W_\eta))>\eta\}\subset \o'$, where $W_\eta:=\{x\in \O|g(x)>\eta\}$
\item[(ii)] for all positive continuous solution $u$ of \eqref{pev.eq.har}, the following inequality holds: 
$$ u(x)\le C u(y) \quad \text{ for all }\quad x \in \o ,y\in \o'\cap \o .$$
\end{itemize}
\end{theorem}

Next, we present a contraction Lemma which guarantees that when $\O$ is bounded then any continuous positive solution $u$ of equation \eqref{pev.eq.har} is bounded in $ \bar \O$.  
\begin{lemma}[Contraction Lemma]\label{pev-lem-esti1}
Let $\O\subset\R^n$ and $u\in C(\O)$ be respectively an open set and a positive solution of \eqref{pev.eq.har}. Then 
there exists $\eps^*>0$ such that for all $\eps\le \eps^*$, there exists $\O_\eps$ and $C(\alpha,\beta,J,\eps,b)$ such that
$$\int_{\O_{\eps}}u(y)\,dy \ge C\int_{\O}u(y)\,dy.$$
Moreover, $\O_\eps$ satisfies the following chain of inclusion 
$$\left\{x\in\O| d(x,\partial\O)>\alpha\eps \right\}\subset\O_\eps\subset\left\{x\in\O| d(x,\partial\O)>\frac{\alpha\eps}{2}\right\}.$$
\end{lemma}
A proof of these  results  can be found in \cite{Co4}.
\medskip

\section{Construction of a principal eigenpair  \label{pev.s.ev}}
In this section we prove the  criterion of existence of a principal eigenpair (Theorems \ref{pev.th1},\ref{pev.th5} and \ref{pev.th7}). That is, we prove the existence of a solution $(\lambda_p,\phi_p)$  of the equation 
\begin{equation}
\label{pev.eq.eigenfunction} \oplb{\phi_p}{\O}+a(x)\phi_p=-\lambda_p\phi_p
\quad \hbox{in $\O$}.
\end{equation}
with $\phi_p>0$, $\phi_p\in C(\O)$ and $\lambda_p$ is the principal eigenvalue of $\lb{\O}+a(x)$ defined by \eqref{pev.eq.def}.  In this task,  we first  restrict our analysis to the case of a bounded domain $\O$ and then prove the criterion for unbounded domains.  We split this section into two subsections, each of them dedicated to one situation. 
\smallskip

\subsection{Existence of a principal eigenpair when $\O$ is a   bounded domain \label{pev-ss.bd}}~

To simplify  the presentation, we will first concentrate our attention on the construction of a principal eigenpair when $J,g,b$ satisfy the assumptions $(H2-H4)$ (Theorem \ref{pev.th1}). Then we provide  an argumentation for  the construction of  a principal eigenpair when the assumptions on $g$ are relaxed  (Theorem \ref{pev.th5}).

In a first step, let us  show that   the eigenvalue problem \eqref{pev.eq.eigenfunction} admits a positive solution i.e. there exists $(\mu_{1,0},\phi_1)$ with $\phi_1>0$, $\phi_1\in L^{\infty}(\O)\cap C(\bar \O)$ solution of \eqref{pev.eq.eigenfunction}. More precisely, we  prove   
\begin{theorem}\label{pev.th.pev}
Let $\O\subset  \R^n$ be a bounded open set and assume that $J,g,$ and  $a(x)$ satisfy $(H1-H4)$. Let us denote $\sigma:=\sup_{\bar \O} a(x)$ and $\O_\theta:=\{x\in\O|d(x,\partial \O)>\theta\}$.
Assume further that the function   $a(x)$ satisfies $\frac{1}{\sigma -a(x)}\not \in L^1_{d\mu}(\bar \O)$. Then there exists $\theta_0>0$ such that for all $\theta\le \theta_0$ the operator $\lb{\O_\theta} + a(x)$ has a unique eigenvalue $\mu_{1,\theta}$ in
$C(\O_\theta)$, that is to say, there is an unique $\mu_{1,\theta} \in \R$  such that 
\begin{equation}
\oplb{\phi_1}{\O_\theta}+a(x)\phi_1=-\mu_{1,\theta}\phi_1 \quad \hbox{in $\O_\theta$}.
\end{equation}
admits a positive solution $\phi_1 \in C(\bar\O_\theta)$. Moreover
$\mu_{1,\theta}$ is simple (i.e the space of $C(\bar\O_\theta)$ solutions to \eqref{pev.eq.eigenfunction} is one dimensional) and satisfies
$$\mu_{1,\theta} <-\max_{\bar\O_\theta} a(x).$$
\end{theorem}

Suppose for the moment that the above Theorem  holds true. To conclude the proof of Theorem \ref{pev.th1} which establishes the criterion of existence of an eigenpair,   we are left to show that  the principal eigenvalue defined by \ref{pev.eq.def} is the same as the one obtained in the Theorem \ref{pev.th.pev} for $\theta=0$. 
Namely, we are  reduced to prove of the following results. 

\begin{lemma}\label{pev.lem.eig-equality}
Let  $a(x)$  be as in Theorem \ref{pev.th.pev} then we have $\lambda_p=\mu_{1,0}$ where $\lambda_p$ and  $\mu_{1,0}$ are respectively the principal eigenvalue of $\lb{\O}+a(x)$ defined by \eqref{pev.eq.def} and the  eigenvalue of $\lb{\O}+a(x)$  obtained in Theorem \ref{pev.th.pev}. 
\end{lemma}

Before proving Theorem \ref{pev.th.pev}, let us prove  the above Lemma.  

\dem{Proof of Lemma \ref{pev.lem.eig-equality} :}
First, let us  define the following quantity
$$\lambda'_p:=\sup\{\lambda\in \R| \, \exists \, \phi >0,\phi \in C(\bar \O) \text{ so that } \oplb{\phi}{\bar\O} +a(x)\phi +\lambda\phi \le 0 \text{ in } \bar\O \}.$$
Obviously $\lambda'_p$ is well defined and is sharing the same properties than $\lambda_p$. Moreover,   we have $\lambda_p'\le\lambda_p$.
Let us now show that $\lambda_p'=\mu_{1,0}$.
First by definition of $\lambda'_p$ we easily have  $\lambda'_p \ge \mu_{1,0}$.   
Now to obtain the equality  $\lambda_p=\mu_{1,0}$ we  argue by contradiction. Assume that $\lambda'_p>\mu_{1,0}$.
By definition of $\lambda'_p$ there exists $\psi>0$, $\psi\in C(\bar\O)$ such that 
\begin{equation}\label{sigma-supersol}
\oplb{\psi}{\O}+(a(x)+\lambda)\psi\le 0 \quad \text{ in }\quad \bar \O.
\end{equation}
Observe that we can rewrite $\oplb{\phi_1}{\O}+a(x)\phi_1$ the following way
\begin{align*}
\oplb{\phi_1}{\O}+a(x)\phi_1 &= \int_{\O}J\left[\frac{x-y}{g(y)}\right]\frac{\phi_1(y)}{g(y)}\, dy+  a(x)\phi_1 \\
 &= \int_{\O}J\left[\frac{x-y}{g(y)}\right]\frac{\psi(y)\phi_1(y)}{\psi(y)g(y)}\, dy+  a(x)\frac{\phi_1(x)}{\psi (x)}\psi(x) 
\end{align*}
From \eqref{sigma-supersol}, we  find that
$$a(x)\psi \le  - \oplb{\psi}{\O} -\lambda \psi$$ 
and it follows that 
 $$
\oplb{\phi_1}{\O}+a(x)\phi_1 \le   \int_{\O}J\left[\frac{x-y}{g(y)}\right]\frac{\psi(y)}{g(y)}\left[\frac{\phi_1(y)}{\psi(y)}-\frac{\phi_1(x)}{\psi(x)}\right]\, dy-  \lambda\frac{\phi_1(x)}{\psi (x)}\psi(x). 
$$
By using the definition of $\mu_{1,0}$, we end up with the following inequality
\begin{equation}
\int_{\O}J\left[\frac{x-y}{g(y)}\right]\frac{\psi(y)}{g(y)}\left[\frac{\phi_1(y)}{\psi(y)}-\frac{\phi_1(x)}{\psi(x)}\right]\, dy \ge (\lambda -\mu_{1,0})\phi_1>0. \label{pev.eq.inequa} 
\end{equation} 
Let us denote  $w:=\frac{\phi_1}{\psi}$. Observe that by \eqref{sigma-supersol} $w  \in L^{\infty}\cap C(\bar\O)$, therefore $w$ achieves a global maximum somewhere in $\bar \O$, say at $\bar x$.
By using the inequality \eqref{pev.eq.inequa} at the point $\bar x$,  we find the following contradiction 
$$    0<\int_{\O}J\left[\frac{\bar x-y}{g(y)}\right]\frac{\psi(y)}{g(y)}[w(y)-w(\bar x)]\, dy\le 0.    $$
Thus $\mu_{1,0}=\lambda'_p$.

Observe now that if there exists a positive eigenfunction $\psi \in C(\O)\cap L^{\infty}(\O)$ associated to the principal eigenvalue $\lambda_p$, i. e. $\oplb{\psi}{\O}+(a(x)+\lambda_p)\psi=0$, then  we have $\psi \in C(\bar \O)$. Therefore, using the definition of $\lambda'_p$  it follows that $\lambda_p\le \lambda'_p=\mu_{1,0}\le\lambda_p$.
To conclude the proof,  we are  left to show that such bounded function $\psi$ exists.

So let $(\theta_n)_{n\in\N}$ be a positive sequence which converges to 0 and consider the sequence  of set $(\O_{\theta_n})_{n\in \N}$ defined in Theorem \ref{pev.th.pev}. 
By construction, using the monotonicity property of the principal eigenvalue with respect to the domain ( (i) of Proposition \ref{pev.prop1}) we deduce that $(\lambda'_p(\lb{\O_{\theta_n}} +a(x)))_{n\in\N} $ is a non increasing bounded sequence. Namely,  we have for all $n\in \N$ $$\lambda_p(\lb{\O} +a(x))\le\lambda'_p(\lb{\O_{\theta_{n+1}}} +a(x)) \le \lambda'_p(\lb{\O_{\theta_n}} +a(x)). $$ 
Thus, as $n$ goes to infinity $\lambda'_p(\lb{\O_{\theta_n}} +a(x))$ converges to some $\bar \lambda \ge \lambda_p$.

On another hand since $\theta_n$ tends to $0$, by Theorem \ref{pev.th.pev}, there exists $n_0$ so that for all $n\ge n_0$, a principal eigenpair $(\mu_{1,\theta_n},\phi_n)$ exists for the operator $\lb{\O_{\theta_n}}+a(x)$.  
Arguing as above, we conclude that  $\mu_{1,\theta_n}=\lambda'_p(\lb{\O_{\theta_n}}+a(x))$.

We claim that
\begin{cla}\label{pev.cla.estin} There exists $n_1 \in \N$ such that for all $n\ge n_1$ we have $\mu_{1,\theta_n}<-\sigma = -sup_{\O}a(x)$. 
\end{cla} 
Assume for the moment that the claim holds. Then the final argumentation goes as follows.
Next, let us normalized $\phi_n$ so that $\sup_{\O_{\theta_n}}\phi_n =1 $.  
With this normalisation  $(\phi_n)_{n\in\N}$ is an uniformly bounded sequence of continuous functions. So by a standard diagonal extraction argument, there exists a subsequence still denoted $(\phi_n)_{n\in \N}$ such that $(\phi_{n})_{n\in\N}$ converges  locally uniformly  to  a non negative  bounded continuous function  $\psi$. Furthermore,  $\psi$   satisfies 

$$ \oplb{\psi}{\O}+(a(x)+\bar \lambda) \psi =0.$$
 
Now recall that  $(\mu_{1,\theta_n},\phi_n)$ satisfies
$$ \oplb{\phi_n}{\O_{\theta_n}}+a(x) \phi_n +\mu_{1,\theta_n}\phi_n=0.$$
Using the above claim, we have  $\mu_{1,\theta_n}<-\sigma = -sup_{\O}a(x)\le-sup_{\O_{\theta_n}}a(x) $ for $n$ big enough, so $\sup_{\O_ {\theta_n}}(a(x)+\mu_{1,\theta_n})<0$ and the uniform estimates i.e. Theorem \ref{pev.th.uniformesti} applies to $\phi_n$. Thus we have for $\eta>0$ small fixed independently of $n$ 
$$1\le C(\eta)\phi_n(x) \quad \text{ for all }\quad  x \in \{x\in \O_{\theta_n}| d(x,\partial \O_{\theta_n})>\eta\}.$$ 
Therefore $\psi$ is   non trivial and $(\bar \lambda,\psi)$ solves the eigenvalue  problem \eqref{pev.eq.eigenfunction}.
 Using once again the equation satisfied by $\psi$ and the definition of $\lambda_p$, we easily obtain that $\bar\lambda \le\lambda_p\le \bar\lambda$ which proves that $\psi$ is our desired eigenfunction associated to $\lambda_p$.
 \fdem

Let us turn our attention to proof of the Claim \ref{pev.cla.estin}. But before proving the Claim let us establish the following a useful estimate.

\begin{lemma}\label{pev.lem.estiJ}
There exists positives constants  $r$ and $c_0$   so that $$\forall x\in \bar \O ,\;  
\int_{B_r(x)\cap\bar \O}J\left(\frac{x-y}{g(y)}\right)u(y)d\mu (y)\ge c_0\int_{B_r(x)\cap\bar \O}u(y)d\mu (y).$$
\end{lemma}

\dem{Proof:}
Since $J$ is continuous and $J(0)>0$, there exists $\delta>0$ and $c_0>0$ so that for all $z\in B(0,\delta)$ we have
$J(z)\ge c_0$.  

 Observe that  for all $(x,y)\in \bar \O\times B_r(x)$  with  $r<\frac{\delta\alpha}{2}$, using that $g\ge \alpha>0$,   we have
 $$  \left\| \frac{x-y}{g(y)}\right\|\le \frac{2r}{\alpha} \le \delta.$$
 Thus, for $r<\frac{\delta\alpha}{2}$  and $y\in B_r(x)$  we have $ J\left(\frac{x-y}{g(y)}\right) >c_0,$ 
and the estimate follows.
 \fdem

We are now in position to prove Claim \ref{pev.cla.estin}.
\dem{Proof of Claim \ref{pev.cla.estin}}
Let us denote $\sigma$  the maximum of $a(x)$ in $\bar \O$. By assumption, we have  $\frac{1}{\sigma -a(x)}\not\in L^1_{loc}(\bar\O)$. So there exists $x_0\in \bar\O$ such that $\frac{1}{\sigma -a(x)}\not\in L^1(B_r(x_0)\cap\bar\O)$ and for $\eps$ small enough say $\eps\le \eps_0$ we have  $$c_0\int_{\bar\O\cap B(x_0,r)}\frac{dx}{-(a(x) -\sigma+\eps)}\ge 4. $$ 

Choose $n_1$ big enough, so that for all $n\ge n_1$, $B_r(x_0)\cap\bar\O_{\theta_n}\not= \emptyset$. For $\eps\le \eps_0$, since $\O_{\theta_n}\to \O$,  we can increase $n_1$ if necessary to achieve for  all $n\ge n_1$ 
 \begin{equation}\label{pev.eq.esti-thetan}
 c_0\int_{\bar\O_{\theta_n}\cap B(x_0,r)}\frac{dx}{-(a(x) -\sigma-\eps)}\ge 2.
\end{equation}
Recall now  that for $n$ big enough, say $n\ge n_2$,  there exists $(\mu_{1,\theta_n},\phi_n)$ that satisfies the equation
$$ \oplb{\phi_n}{\O_{\theta_n}}+a(x) \phi_n +\mu_{1,\theta_n}\phi_n=0.$$

Since $\phi_n$ is positive we have 
$$  \oplb{\phi_n}{\bar\O_{\theta_n}\cap B(x_0,r)}\le -(a(x) +\mu_{1,\theta_n})\phi_n.$$
Using the Lemma \ref{pev.lem.estiJ}, we see that 
$$ \frac{c_0}{-(a(x) +\mu_{1,\theta_n})}\int_{\bar\O_{\theta_n}\cap B(x_0,r)}\phi_n(y)\,dy\le \phi_n(x).$$
Integrating the above inequality on $\bar\O_{\theta_n}\cap B(x_0,r)$ it follows that 
\begin{align*}
 \int_{\bar\O_{\theta_n}\cap B(x_0,r)}\left(\frac{c_0}{-(a(x) +\mu_{1,\theta_n})}\int_{\bar\O_{\theta_n}\cap B(x_0,r)}\phi_n(y)\,dy\right)&\le \int_{\bar\O_{\theta_n}\cap B(x_0,r)}\phi_n(x)
 \\
\int_{\bar\O_{\theta_n}\cap B(x_0,r)}\left(\frac{c_0}{-(a(x) +\mu_{1,\theta_n})}\right)\int_{\bar\O_{\theta_n}\cap B(x_0,r)}\phi_n(y)\,dy&\le \int_{\bar\O_{\theta_n}\cap B(x_0,r)}\phi_n(x).
\end{align*}
Thus, $$ \int_{\bar\O_{\theta_n}\cap B(x_0,r)}\left(\frac{c_0}{-(a(x) +\mu_{1,\theta_n})}\right)\le 1$$
From \eqref{pev.eq.esti-thetan}, it follows that for all $n\ge \sup(n_1,n_2)$ we have 
$$\mu_{1,\theta_n}\le -\sigma-\eps.$$
\fdem

\begin{remark}\label{pev.rem.esti}
Observe that if $\sup_{\O}a(x)$ is achieved in $\O$ then the estimation $\mu(1,\theta)$ follows immediately from the monotonicity properties of the principal eigenvalue. Indeed,  for $\theta$ small enough, say $\theta \le \theta_0$ we have $\sup_{\O_\theta}a(x)=\sup_ {\O}a(x).$
Hence, $$ \lambda'_p(\lb{\O_\theta}+a(x))\le\lambda'_p(\lb{\O_{\theta_0}}+a(x))< - \sup_{\O_{\theta_0}}a(x)=-\sigma.$$
\end{remark}

Let us now turn our attention to the proof of Theorem \ref{pev.th.pev}.

 For convenience, in this proof we write the eigenvalue problem
\begin{align*}
\oplb{u}{\O_\theta}+a(x)u  = -\mu u
\end{align*}
in the form
\begin{align}
\label{pev-eq.pev2} \oplb{u}{\O_\theta} + \bar a(x) u = \rho u
\end{align}
where
$$  \bar a(x) = a(x) + k, \quad \rho = -\mu  + k$$
and $k>0$ is a constant such that $\inf_{\O_\theta} \bar a >0$.

Let us now prove the following  useful result:
\begin{lemma}\label{pev.lem.tec}
Let $\O,J,g$ and $a$ be as in Theorem \ref{pev.th.pev}. 
Then there exists $\theta_0>0$ so that for all $\theta\le \theta_0$ there exists $\delta>0$
and $u \in C(\bar\O_\theta)$, $u \ge 0$, $u \not \equiv 0$, such that 
$$
\lb{\O_\theta}[u]+\bar a(x)u\ge (\bar \sigma+\delta) u,
$$
where $\bar \sigma(\theta):=\max_{ \bar \O_\theta} \bar a(x)$.
\end{lemma}

Observe that the proof of  Theorem \ref{pev.th.pev} 
easily follows from the above Lemma.  Indeed, if the Lemma holds true,
since under the assumption (H1-H4) the operator  $\lb{\O}:C(\bar\O_\theta) \to C(\bar\O_\theta)$ is  compact,   we have
$r_e(\lb{\O_\theta}+\bar a(x))=r_e(\bar a(x)) {=} \bar \sigma(\theta)$. Thus
$(\bar \sigma(\theta)+\delta)>r_e(\lb{\O_\theta}+\bar a(x))$ and the existence  Theorem of Edmund \textit{ et al.} (Theorem
\ref{pev.th.eps}) applies.

 Finally we observe that the principal eigenvalue is simple
since for a bounded domain $\O$ the cone of positive continuous functions has a non-empty
interior and, for a sufficiently large $p$, the operator $(\lb{\O_\theta} +
\bar a)^p $ is strongly positive, that is, it maps $u \ge 0$, $u
\not\equiv0$ to a strictly positive function, see \cite{zeidler}. 
\fdem

\begin{remark}
Note that the simplicity of the eigenvalue $\mu_{\theta}$ requires that $\O_\theta$ is a connected set. Indeed, when open set $\O$ is not connected,  it may happen that the operator $(\lb{\O_\theta} +\bar a)^p $ is never strongly positive in $C(\bar\O)$ and several non-positive eigenfunction exists with no positive eigenfunction. 
\end{remark}
%
%
%

Let us now turn our attention to the proof of  Lemma \ref{pev.lem.tec}: 
\dem{Proof of the Lemma \ref{pev.lem.tec}:} 
Let  us denote $\Gamma$ the closed set  where the continuous function $\bar a$ takes its maximum $\bar\sigma$ in $\bar \O$.
$$\Gamma:=\{z\in \bar \O |\; \bar a(z)=\bar \sigma\}.$$
Since $\bar a$ is a continuous function and $\O$ is bounded, $\Gamma$ is a compact set. Therefore $\Gamma$ can be covered by a finite number of balls of radius $r$, i.e. $ \Gamma \subset \bigcup_{i=1}^{N}B_r(x_i)$  with $x_i\in \G $.
By construction, we have   $\frac{1}{\bar \sigma -\bar a(x)}=\frac{1}{\sigma -a(x)} \not\in L^{1}_{d\mu,loc}(\bar \O)$. Therefore  
$\frac{1}{\bar\sigma -\bar a(x)} \not \in L^1_{d\mu}(\bigcup_{i=1}^{N}B_r(x_i)\cap \bar \O)$ and there exists $-\lambda_0>\bar\sigma$   so that  for some $x_i$ we have
\begin{equation}
\int_{B_r(x_i)\cap\bar  \O}\frac{c_0 }{-\lambda_0 -\bar a(x)}d\mu\ge 4.
\end{equation}
Since $\O_\theta \to \O$ as $\theta$ tends to 0 there exists $\theta_0$ so that for all $\theta\le \theta_0$ we have 
 \begin{equation}
\int_{B_r(x_i)\cap\bar  \O_\theta}\frac{c_0 }{-\lambda_0 -\bar a(x)}d\mu\ge 2.
\label{pev.eq.lambda_0}
\end{equation}
Let us fix $x_i$ such that \eqref{pev.eq.lambda_0} holds true and let us denote $\o_{_\theta}:= B_r(x_i)\cap \bar \O_\theta$.
 We consider now the following eigenvalue problem
\begin{equation}\label{pev.eq.exist}
c_0\int_{\o_{_\theta}}u(y)\,d\mu(y) +\bar a(x)u(x)+\lambda u(x)=0,
\end{equation}
where $c_0$ is the constant obtained in the Lemma \ref{pev.lem.estiJ}.

We claim that 
\begin{cla}\label{pev.cla.exist}
There exists $(\lambda_1, \phi_1)$  solution of \eqref{pev.eq.exist} so that $\phi_1\in L^{\infty}(\o_{_\theta})\cap C(\o_{_\theta})$ and $\phi_1>0$.
\end{cla}

Observe that by proving this claim we end the proof of the Lemma. Indeed, fix $\theta <\theta_0$ and assume for the moment that this claim holds true. Then there exists 
$(\lambda_1, \phi_1) $ such that 
\begin{equation}\label{pev.eq.exist1}
c_0\int_{\o_{_\theta}}\phi_1(y)\,d\mu(y) +\bar a(x)\phi_1(x)+\lambda_1 \phi_1(x)=0.
\end{equation}
Obviously, for any positive constant $\rho$, $(\lambda_1,\rho\phi_1)$ is also a solution of the equation \eqref{pev.eq.exist1}. Therefore without any loss of generality we can assume that $\phi_1$  is such that $\phi_1\le 1$.
Set $\tilde c_0:=c_0\int_{\o_{_\theta}}\phi_1(y)d\mu(y)$. From the equation \eqref{pev.eq.exist1},  since $0<\phi_1\le 1$ we see easily  that  
\begin{equation*}
 -(\lambda_1+\bar a(x))>\tilde c_0.
\end{equation*}
Therefore there exists a positive constant $d_0$ such that 
\begin{equation}\label{pev.eq.exi.def1}
\phi_1\ge d_0 \quad \text{ in }\quad \o
\end{equation}
and 
\begin{equation}
 -(\lambda_1+\bar\sigma(\theta))\ge \tilde c_0>0.\label{pev.eq.exi.def2}
\end{equation}

Let us now consider a set $\o_\eps\subset\subset \o_{_\theta}$  which verifies 
\begin{equation}\label{pev.eq.exist2}
\int_{\o_{_\theta}\setminus \o_\eps}\,d\mu\le \frac{d_0|\lambda_1+\bar\sigma(\theta)|}{2c_0}.
\end{equation}

Since by construction $\bar \O_\theta\setminus \o_{_\theta}$ and $\bar \o_\eps$ are two disjoint  closed  subsets of $\O_\theta$,  the  Urysohn's Lemma applies and there exists a positive  continuous  function $\eta$ such that $0\le \eta \le 1$, $\eta(x) = 1$ in $\o_\eps$, $\eta(x) = 0$ in $\bar \O_\theta\setminus \o_{_\theta}$. 

Next, we define  $w:=\phi_1\eta $ and  we compute $\oplb{w}{\O_\theta}+b(x)w$.

Since $w\equiv 0$ in  $\bar\O_\theta\setminus \o_{_\theta}$, 
we have 
$$\oplb{w}{\O_\theta}+\bar a(x)w=\int_{\o_{_\theta}}J\left(\frac{x-y}{g(y)}\right) w(y)\, d\mu \ge (\bar\sigma(\theta)+\delta ) w =0    $$
for any $\delta>0$.

On another hand, in $\o_{_\theta}$,  by using the Lemma \ref{pev.lem.estiJ}  we see that 
\begin{align}
\oplb{w}{\O_\theta}+\bar a(x)w&=\int_{\o_{_\theta}}J\left(\frac{x-y}{g(y)}\right) w(y)d\mu +\bar a(x)w\\ 
&\ge c_0\int_{\o_{_\theta}}w(y)d\mu(y) + \bar a(x)w\\
&\ge c_0\int_{\o_\eps}\phi_1 d\mu(y) + \bar a(x)w.\label{pev.eq.exist3}
\end{align}
 Since $\phi_1$ satisfies the equation \eqref{pev.eq.exist1}, using the estimates \eqref{pev.eq.exi.def1}, \eqref{pev.eq.exi.def2} and  \eqref{pev.eq.exist2} we deduce from the inequality \eqref{pev.eq.exist3}   that      
\begin{align}
\oplb{w}{\O_\theta}+\bar a(x)w&\ge -(\lambda_1 +\bar a(x))\phi_1 +\bar a(x)w -c_0\int_{\o_{_\theta}\setminus \o_\eps}\phi_1\\
& \ge\frac{|\lambda_1+\bar\sigma(\theta)|}{2}\phi_1 +(\bar\sigma(\theta) -\bar a(x))\phi_1 +\bar a(x)w + \frac{d_0|\lambda_1+\bar\sigma(\theta)|}{2}-c_0\int_{\o_{_\theta}\setminus \o_\eps}\phi_1 \, d \mu.\\
&\ge \left(\frac{|\lambda_1+\bar\sigma(\theta)|}{2}\right)\phi_1 +(\bar\sigma(\theta) -\bar a(x))\phi_1 +\bar a(x)w.\label{pev.eq.exist4}
\end{align}
where we use in the last inequality, that $\phi_1\le 1$ and the estimate \eqref{pev.eq.exist2}.

%

Since $(\bar\sigma(\theta)-\bar a(x))$ and $\frac{|\lambda_1+\bar\sigma(\theta)|}{2}$ are two positive quantities and $\phi_1\ge w$,  we conclude that   
 \begin{equation}
\oplb{w}{\O_\theta}+\bar a(x)w\ge \left(\frac{|\lambda_1+\bar\sigma(\theta)|}{2} + \bar\sigma(\theta)\right) w.
\end{equation}

Hence, in $\O_\theta$,  $w$ satisfies $$ \oplb{w}{\O_\theta}+\bar a(x)w\ge (\bar\sigma(\theta)+\delta)w, $$ 
with $\delta =\frac{|\lambda_1+\bar\sigma(\theta)|}{2}$, which proves the Lemma.
\fdem

Let us now prove the claim \ref{pev.cla.exist}.
\dem{Proof of Claim \ref{pev.cla.exist}}
Fix $\theta \le \theta_0$. For $\lambda<-\bar\sigma(\theta)$, consider the positive function $\phi_\lambda:= \frac{c_0}{-\lambda -\bar a(x)}$.  Let us substitute $\phi_\lambda$ into the equation \eqref{pev.eq.exist}, then we have   
$$c_0\int_{\o_{_\theta}} \phi_\lambda\, d\mu -c_0=0.$$
Therefore, we end the proof of the Claim \ref{pev.cla.exist} by finding $\lambda$ such that $\int_{\o_{_\theta}} \phi_\lambda\, d\mu =1$. 
Observe that the functional $F(\lambda):=\int_{\o_{_\theta}} \phi_\lambda\, d\mu$ is continuous and monotone increasing with respect to $\lambda$ in $(-\infty,-\bar\sigma)$. 
Moreover, by construction, we have:  
$$ \lim_{\lambda \to -\infty} F(\lambda)= 0 \quad \text{ and }\quad  F(\lambda_0)\ge 2. $$
Hence by continuity there exists a  $\lambda_1$ such that $F(\lambda_1)=1$.
\fdem

Now we expose the argumentation for the construction of a principal eigenpair when the assumptions on $g$ are relaxed and prove the Theorem \ref{pev.th5}. 
To show Theorem \ref{pev.th5} we follow the scheme of the argument developed above. 
\dem{Proof of Theorem \ref{pev.th5}} 
As above,  we can rewrite the eigenvalue problem \eqref{pev.eq.eigenfunction} the following way  

\begin{align}
\label{pev.eq.pev3} \oplb{u}{\O_\theta} + \bar a(x) u = \rho u
\end{align}
with
$$  \bar a(x) = a(x) + k, \quad \rho = -\mu  + k$$
and $k>0$ is a constant such that $\inf_{\O_\theta} \bar a >0$.

Observe that under the assumptions $(H1, H2,\tilde H3,H4)$  the following  family
$$\lb{\O_\theta}(B_1):=\{ \oplb{f}{\O_\theta} \; / \; f : \O \to \mathbb{R}, \;\; ||f||_\infty \leq 1 \}$$
is equicontinuous.  Indeed, let $\eps >0$ be fixed.  Since  $\frac{1}{g^n}\in L^p_{loc}(\bar\O_\theta)$, there exists $\eta
>0$ such that
\begin{equation}\label{equi-1}
\int_{_{\O_\theta\cap\{g<\eta\}}} \frac{dy}{g^n(y)} <\frac{\eps}{4||J||_\infty}.
\end{equation}

From the uniform continuity of   $J$ in the unit ball $B(0,1)$, we deduce that  there
exists $\gamma>0$ such that for $|w- \bar w|<\gamma/\eta$,
\begin{equation}\label{equi-2} 
|J(w)-J(\bar w)|<\eps\eta^n/2|\O_\theta|.
\end{equation} 
A short computation using \eqref{equi-1} and \eqref{equi-2} shows that  for $|x-z|<\gamma$ 
\begin{align*}
|\oplb{f}{\O_\theta}(x)-\oplb{f}{\O_\theta}(z)|&\le  \int_{\O_\theta}\left | J\left[\frac{x-y}{g(y)}\right] - J\left[\frac{z-y}{g(y)}\right]     \right | \left |\frac{f(y)}{g^n(y)}\right |\,dy\\
&\le  2||J||_\infty\int_{_{\O_\theta\cap\{g<\eta\}}} \frac{1}{g^n(y)}dy +\frac{1}{\delta^n}\int_{_{\O_\theta\cap\{g\ge\eta\}}}\left | J\left[\frac{x-y}{g(y)}\right] - J\left[\frac{z-y}{g(y)}\right]     \right | \,dy\\
&\le \eps.
\end{align*}
Hence, $\lb{\O_\theta}(B_1)$ is equicontinuous and $\lb{\O_\theta}:C(\bar\O_\theta) \to C(\bar\O_\theta)$ is a  compact operator.  

Next, we show  the following 
\begin{lemma}\label{pev.lem.tec2}
Let $\O,J,g$ and $a$ be as in Theorem \ref{pev.th5}. Then there exists $\theta_0$ so that for all $\theta\le\theta_0$
  there exists $\delta>0$
and $u \in C(\bar\O_\theta)$, $u \ge 0$, $u \not \equiv 0$, such that
$$
\lb{\O_ \theta}[u]+\bar a(x)u\ge (\bar \sigma+\delta) u.
$$
\end{lemma}

 As above the existence of a positive eigenpair $(\rho,\phi)$ easily follows from the Lemma \ref{pev.lem.tec2}. 
 Arguing as above, we see that   $\mu_{1,0}=\lambda_p(\lb{\O}+a(x))$, which concludes the proof of Theorem \ref{pev.th5}.  

\fdem

Let us turn our attention to the proof of Lemma \ref{pev.lem.tec2}
\dem{Proof of Lemma \ref{pev.lem.tec2}}
First let us recall that by assumption $\stackrel{\circ}{\G}\neq \emptyset$ where $\G:=\{x\in \bar \O | a(x)=\sigma\}$ and let us  define 
the following set $\Sigma_\eta:=\{x\in \O| g(x)\ge \eta\}$. 

By construction, we easily see that $\stackrel{\circ}{\G'}\neq \emptyset$ where $\G':=\{x\in \bar \O | \bar a(x)=\bar \sigma\}$.
Therefore, there exists $x_0\in \O$ and $\eps>0$ such that $B_{\eps}(x_0)\subset (\stackrel{\circ}{\G'}\cap \O)$.
Moreover for $\theta$ small, say $\theta\le \theta_0$ we have $B_{\eps}(x_0)\subset (\stackrel{\circ}{\G'}\cap \O_\theta)$.

Let us define $\o_\eta:= B_\eps(x_0)\cap \Sigma_\eta$. By assumption we have $\frac{1}{g^n}\in L^{p}(\O)$, so for $\eta$ small enough $\o_\eta$ is a non void open subset of $\O_\theta$ for $\theta\le \theta_0$.

Let us now consider the eigenvalue problem \eqref{pev.eq.pev3} with $\O=\o_\eta$, i.e
$$ \oplb{u}{\o_\eta} + \bar a(x) u = \rho u \quad \text { in }\quad \o_\eta.$$
By construction,  in $B_\eps(x_0)$ we have  $\bar a(x)\equiv \bar\sigma$. So the above equation reduces to:
\begin{align}\label{pev.eq.pev-eps}
   \oplb{u}{\o_\eta} =  \bar \rho u \quad \text { in }\quad \o_\eta,
\end{align}
where $\bar \rho=(\rho-\bar\sigma)$.

Since $\lb{\o_\eta}$ is a compact strictly positive operator in $C(\bar\o_\eta)$, using Krein-Rutmann Theorem there exists a positive  eigenvalue $\bar \rho_1>0$ and a positive eigenfunction $\phi_1 \in C(\bar\o_\eta) $ such that  $(\bar\rho_1, \phi_1)$ satisfies \eqref{pev.eq.pev-eps} i.e
  $$\oplb{\phi_1}{\o_\eta} =  \bar \rho \phi_1.$$

Arguing as in Lemma \ref{pev.lem.tec}, for all $\theta\le \theta_0$ we can construct  a nonnegative test function $u$ such that 
$$  \oplb{u}{\O_\theta}+\bar a(x)u\ge (\delta+\bar\sigma)u,$$
for a $\delta>0$ small enough.  

\fdem

\begin{remark}
Observe that all the previous constructions can be easily adapted to an operator $\t +a(x)$ where $\t$ is an integral operator with a continuous nonnegative kernel $k(x,y)$ that satisfies $\tilde H2$, i.e.
$$
\exists\, c_0>0,\ \eps_0 >0\; \text{ such that } \min_{x\in\O}\left(\min_{y\in B(x,\eps_0)}k(x,y)\right)>c_0.
$$
In particular, we can extend the criterion of existence of a principal eigenpair for an operator $\t +a(x)$ where $\t$ is an integral operator with a  kernel $k(x,y)$  that only satisfies that  there exists a  positive integer $N$, so that the kernel  $k_N(x,y)$  satisfies $(\tilde H2)$ where $k_N$ is defined by the recursion:
 \begin{align*}
&k_1(x,y):=k(x,y)\\
&k_{N+1}(x,y):=\int_{\O} k_{N}(x,z)k_1(z,y)\,dz \quad \text{for } \quad N\ge 1.
\end{align*}
  Indeed, in this situation the construction of a test function $u$ (Lemma \ref{pev.lem.tec} or Lemma \ref{pev.lem.tec2} ) holds also for the operator $ \t^N+\bar a^N(x)$. Using that $\bar a \ge 0$, we deduce   
$$(\t+\bar a(x))^N[u]\ge \t^N u+\bar a^N(x)u\ge (\bar \sigma^N +\delta)u.$$
Since in this situation $\t$ is a compact operator, we also have $r_e((\t+\bar a(x))^N)=r_e(\bar a(x)^N)$. Thus
$(\bar \sigma^N+\delta)>r_e((\t+\bar a(x))^N)$ and the Theorem
\ref{pev.th.eps} applies.  Hence, there exists an unique principal eigenpair $(\lambda_p,\phi_p)$ of the following problem
$$ (\t+\bar a(x))^N\phi_p=-\lambda_p\phi_p$$

To obtain a principal eigenpair for $\t +a$ we argue as follows. 
Applying $\t +a(x)$ to the above equation it follows that 
\begin{align*}
(\t+\bar a(x))^{N+1}\phi_p&=-\lambda_p(\t +a(x))\phi_p\\
(\t+\bar a(x))^{N}\psi&=-\lambda_p\psi
\end{align*}
with $\psi:=(\t +a(x))\phi_p$. Since $(\t+\bar a)^N$ is positive operator in $C(\bar\O)$, $\lambda_p$ is simple, we have 
$\psi=\rho\phi_p$. Hence,   $((-\lambda_p)^{\frac{1}{N}},\phi_p)$ is the  principal eigenpair of $\t+\bar a(x)$. 
\end{remark}
\medskip

\subsection{Construction of a Principal eigenpair when $\O$  is an unbounded domain}~

For simplicity in the presentation of the arguments and since the proof of the existence of a principal eigenpair  under the relaxed assumptions  does not significantly differ, we will only present the case where  $\O,J,g$ and $a$ satisfy the assumptions $(H1-H4)$.

To construct  an eigenpair $(\lambda_p,\phi_p)$ in this situation, we proceed using a standard approximation scheme .
 
First let us recall that,  by assumption,  there exists $\O_0\subset \bar \O$ a bounded subset such that $\frac{1}{\sigma -a(x)}\not\in L^1_{d\mu}(\bar  \O_0)$. 
Let $(\o_n)_{n\in \N}$ be a sequence of  bounded increasing  set which covers $\O$, i.e.

$$\o_n \subset \o_{n+1},\quad \bigcup_{n\in \N} \o_n=\O.$$ 
Without loss of generality, we can also assume that $\O_0\subset\o_0$ and therefore  $\frac{1}{\sigma -a(x)}\not\in L^1_{d\mu}(\bar  \o_n)$ for all $n\in \N$. 
 Observe that for each  $\o_n$ the Theorem \ref{pev.th.pev} and the Lemma \ref{pev.lem.eig-equality} apply.  Therefore for each $n$ there exists  a principal eigenpair $(\lambda_{p,n},\phi_{p,n})$ to the eigenvalue problem \eqref{pev.eq.eigenfunction} with $\o_n$ instead of $\O$.

By construction, using the monotonicity of the sequence of $(\o_n)_{n\in\N}$ and the assertion (i) of the Proposition \ref{pev.prop1} we deduce that 
$(\lambda_{p,n})_{n\in\N}$ is a monotone non increasing sequence which is bounded from below. Thus 
$\lambda_{p,n}$ converges to some $\bar\lambda\ge\lambda_p(\lb{\O}+a(x)).$ 
 Moreover, we also have that for all $n\in\N$ 
$$\lambda_p(\lb{\O}+a(x))\le \bar\lambda \le \lambda_{p,n}<\lambda_{p,0}<-\sup_{\bar \O}a(x)=\sigma.$$  

Let us  now fix  $x_1\in \o_0\cap \O$. Observe that since for each integer $n$ the eigenvalue $\lambda_{p,n}$ is simple we can  normalize  $\phi_{p,n}$ by $\phi_{p,n}(x_1)=1$. 

Let us now define   $b_n(x):=-\lambda_{p,n}-a(x)$. Then $\phi_{p,n}$ satisfies  
\begin{equation}\label{pev.eq.approx1}
  \oplb{\phi_{p,n}}{\o_n}=b_n(x)\phi_{p,n} \quad {in}\quad \o_n.
  \end{equation}
 By construction for all $n\in \N$ we have $b_n(x)\le-\lambda_{p,0}-\sigma>0 $,  therefore the Harnack inequality (Theorem \ref{pev.th.harnack}) applies to $\phi_{p,n}$. Thus for $n$ fixed and for all compact set $\o' \subset \subset \o_n$ there exists a constant $C_n(\o')$ such that 
$$\phi_{p,n}(x)\le C_n(\o')\phi_{p,n}(y) \quad \forall \quad x,y \in \o'.$$

Moreover the constant $C_n(\o')$ only depends  on $\bigcup_{x\in \o}B(x,\beta)$ and is  monotone decreasing  with respect to  $\inf_{x\in \o_n}b_n(x)$.
For all $n$, the function $b_n(x)$ being uniformly bounded from below by a constant independent of n, the constant $C_n$ is bounded from above independently of $n$ by a constant $C(\o')$.  Thus we have 
$$\phi_{p,n}(x)\le C(\o')\phi_{p,n}(y) \quad \forall \quad x,y \in \o'.$$ 

From a standard argumentation, using the normalization $\phi_{p,n}(x_1)=1$, we deduce that  the sequence $(\phi_{p,n})_{n\in\N}$  is bounded in $C_{loc}(\O)$ topology. 
Moreover, from a standard diagonal extraction argument, there exists a subsequence still denoted $(\phi_{p,n})_{n\in \N}$ such that $(\phi_{p,n})_{n\in\N}$ converges  locally uniformly  to  a continuous function  $\phi$. Furthermore,  $\phi$  is a nonnegative non trivial function and $\phi(x_1)=1$.    

Since $J$ has a compact support we can pass to the limit in the equation \eqref{pev.eq.approx1} using the Lebesgue  monotone convergence theorem   and  get 
$$ \int_{\O} J\left(\frac{x-y}{g(y)}\right)\phi(y)d\mu(y)+(\bar\lambda+a(x))\phi(x) =0 \quad{ in }\quad \O.$$
 
 As above using the equation, we deduce that $\phi>0$ in $\O$.  Lastly, from the definition of $\lambda_p$ using $(\bar \lambda,\phi)$ as a test function, we see that $\bar \lambda\le \lambda_p \le \bar \lambda$. 
Hence, $(\bar\lambda,\phi)$ is our desired eigenpair.

\fdem

\begin{remark}
Note that our proof of the existence of a principal eigenpair in this situation relies only on the Harnack estimate which for some form holds true when the assumption on $J$ and $g$ are relaxed.  
\end{remark}
\begin{remark}
From the above proofs, using the properties of the principal eigenvalue,  we can derive a practical   dichotomy for $\lambda_p$. Indeed, either $\lambda_p=-\sigma$ or $\lambda_p<-\sigma$ and  there exists a principal positive eigenfunction $\phi_p$ associated to $\lambda_p$.
\end{remark}

\section{Existence of a Maximum principle  \label{pev.s.mp}}

In this section, we explore the relation between a maximum principle property satisfied by an operator $\m$ and the sign of its principal eigenvalue.  Namely, we prove the  Theorem \ref{pev.th2} that we recall below

\begin{theorem}\label{pev.th.mp}
Assume that  $\O$ is a bounded set and let $J$, $g$ and $a$ be as in Theorem \ref{pev.th1}. Then the maximum principle is satisfied  by $\mb{\O}$ if and only if  $\lambda_p(\mb{\O})\ge 0$. 
\end{theorem}

\dem{Proof of Theorem \ref{pev.th.mp}}
Assume first that the operator satisfies the maximum principle.  From the Theorem \ref{pev.th1}, there exists $(\lambda_p,\phi_p)$ such that $\phi_p\in C(\bar\O)$, $\phi_p>0$ and
 $$\oplb{\phi_p}{\O} +a(x)\phi_p +\lambda_p\phi_p=0.$$ 
 As in the previous section, we have can normalise $\phi_p$ so that we have $1\ge \phi_p\le c_0$.  Furthermore, there exists $\delta>0$ so that $-\lambda_p-\sigma\ge \delta>0$ where $\sigma$ denotes the maximum of $a$ in $\bar\O$.

Assume by contradiction that $\lambda_p< 0$ we have 
$$\oplb{\phi_p}{\O} +a(x)\phi_p = -\lambda_p\phi_p> 0.$$ 

Let us choose $\o\subset \subset \O$ such that $$\int_{\O\setminus \o}d\mu(y)\le\frac{c_0\inf\{\delta,\lambda_p\}}{2\|J\|_\infty}.$$
As in the previous section, we can construct a continuous function $\eta$ such that $0\le \eta \le 1$, $\eta(x) = 1$ in $\o$, $\eta(x) = 0$ in $\partial \O$.
Consider now $\phi_p\eta$ and let us compute $\oplb{\phi_p\eta}{\O}+a(x)\phi_p\eta.$ Then we have 
\begin{align*}
\oplb{\phi_p\eta}{\O}+a(x)\phi_p\eta &\ge-\lambda_p\phi_p -\|J\|\int_{\O\setminus \o}d\mu(y)-a(x)\phi_p(1-\eta)\\
&\ge-\lambda_p\phi_p -\frac{c_0\inf\{\delta,\lambda_p\}}{2}-a(x)\phi_p(1-\eta)\\
&\ge-\lambda_p\phi_p -\frac{c_0\inf\{\delta,\lambda_p\}}{2}-\max\{\sigma,0\}\phi_p\\
&\ge-(\lambda_p+\max\{\sigma,0\})\phi_p -\frac{c_0\inf\{\delta,\lambda_p\}}{2}.
\end{align*} 
Since by assumption $-\lambda_p>0$ and $-\lambda_p-\sigma\ge 0$ it follows from the above inequality that
\begin{align*}
\oplb{\phi_p\eta}{\O}+a(x)\phi_p\eta &\ge-(\lambda_p+\max\{\sigma,0\})c_0 -\frac{c_0\inf\{\delta,\lambda_p\}}{2}\\
&\ge \frac{c_0\inf\{\delta,\lambda_p\}}{2}\ge 0.
\end{align*}
By construction we have $\phi_p\eta \in C(\O)$ that satisfies  
\begin{align*}
&\oplb{\phi_p\eta}{\O} +a(x)\phi_p\eta \ge 0 \quad\text{ in } \quad  \O \\
&\phi_p\eta= 0 \quad\text{ on } \quad \partial \O
\end{align*}
Therefore, by the maximum principle \ref{pev.def.pm}, $\phi_p\eta\le 0$ in $\O$ which is a contradiction. 
Hence, $\lambda_p \ge 0$. 
 
Let us now show the converse implication. Assume that   $\lambda_p(\lb{\O} +a(x))\ge 0$, then we will show that the operator satisfies the maximum principle.
Let $u\not\equiv 0$, $u\in C(\bar\O)$ such that $u\ge 0$ on $\partial \O$ and
$$\oplb{u}{\O}+a(x)u\le 0.$$ Let us show that $u>0$ in $\O$.

By Theorem \ref{pev.th1}, there exists $\phi_p>0$ such that 
$$\oplb{\phi_p}{\O}+a(x)\phi_p=-\lambda_p\phi_p\le 0.$$
Let us rewrite $\oplb{u}{\O} +a(x)u$ the following way
\begin{align*}
\oplb{u}{\O} +a(x)u&=\int_{\O}J\left[\frac{x-y}{g(y)}\right]\frac{\phi_p(y)}{g(y)}\frac{u(y)}{\phi_p(y)}\, dy  + a(x)\phi_p(x)\frac{u(x)}{\phi_p(x)}\\
&=\int_{\O}J\left[\frac{x-y}{g(y)}\right]\frac{\phi_p(y)}{g^n(y)}\left(\frac{u(y)}{\phi_p(y)}-\frac{u(x)}{\phi_p(x)}\right)\, dy -\lambda_p \phi_p \frac{u(x)}{\phi_p(x)}
\end{align*}
Let us set $w:=\frac{u}{\phi_p}$, then we have the following inequality  in $\O$
$$ \int_{\O}J\left[\frac{x-y}{g(y)}\right]\frac{\phi_p(y)}{g^n(y)}\left(w(y)-w(x)\right)\, dy -\lambda_p \phi_p w(x)\le 0.$$
From the above inequality we deduce that  $w$ cannot achieve a non positive  minimum in $\O$ without being constant.
Therefore it follows that either $w>0$ in $\O$ or $w \equiv 0$. Since $u\not\equiv 0$, we have $w>0$.
Hence, $\frac{u}{\phi_p}>0$ which implies that $u>0$. 
\fdem

\begin{remark}
From the proof, we can observe that to show the implication 
\begin{center}"$\lambda_p(\lb{\O}+a(x))>0 \Longrightarrow\lb{\O} +a(x)$ satisfies the maximum principle" 
\end{center}
we do not need  the existence of a principal eigenfunction $\phi_p$ when $\lambda_p(\lb{\O}+a(x))>0$. Indeed, in this situation   we can replace in our argumentation the principal eigenfunction $\phi_p$ by  a well chosen positive function $\psi$  i.e. $\psi>0$ such that   there exists  $0<\lambda\le \lambda_p$ satisfying $\oplb{\psi}{\O}+(a(x)+\lambda)\psi \le 0$ which is always  possible  since  $\lambda_p(\lb{\O}+a(x))>0$. 
 \end{remark}

\section{A counter example \label{pev.s.ce}}
In this section, we provide an example of nonlocal equation where no positive bounded eigenfunction exists. Let  $\O$ be a bounded domain and let us  consider the following  principal eigenvalue problem:  
\begin{align}
\rho\int_{\O}u\,dx + a(x)u=\lambda u, \label{pev-eq-ce}
\end{align}
 where $\sigma =a(x_0)=\max_{\bar \O} a(x)$, $\rho$   is  a positive constant and $a(x)\in C^{0}(\bar\O)$ satisfies the condition $\frac{1}{\sigma -a(x)}\in L^1_{loc}(\O)$.
For this eigenvalue problem, we show the following result
\begin{theorem}
If $\rho$ is so that $\rho\int_{\O}\frac{1}{\sigma -a(x)}<1$, then there exists no bounded continuous positive  principal eigenfunction $\phi$ to \eqref{pev-eq-ce}.
 \end{theorem}

\dem{Proof :}
We argue by contradiction. Let us assume that there exists a bounded positive continuous  eigenfunction $\phi$ associated with $\lambda_p$ that we normalize  by $\int_{\O}\phi =1$. 
By substituting $\phi$ into the equation \eqref{pev-eq-ce} it follows that
$$\rho=(\lambda_p-a(x))\phi.$$
Since $\rho>0$, from the above equation we conclude that $\lambda_p-\sigma \ge\tau>0$.
Therefore $$\phi=\frac{\rho}{\lambda_p-a(x)}.$$
 Next, using  the normalization we obtain
$$1=\rho\int_{ \O}\frac{dx}{\lambda_p-a(x)}.$$
By construction $\lambda_p\ge \sigma$, therefore we have $$1=\rho\int_{\O}\frac{dx}{\lambda_p-a(x)}\le  \rho\int_{\O}\frac{dx}{\sigma-a(x)}.$$
Since $\rho\int_{ \O}\frac{dx}{\sigma-a(x)}<1$ we end up with the following contradiction
  $$1=\rho\int_{\O}\frac{dx}{\lambda_p-a(x)}\le \rho\int_{\O}\frac{dx}{\sigma-a(x)}<1.$$
Hence there exists no positive bounded eigenfunction $\phi$ associated to $\lambda_p$.
\fdem


\section{Existence/Non existence of solution of \eqref{pev.eq.semilin}:\label{pev.s.cri}}
In this section we prove the Theorem \ref{pev.th3} . That is to say, we investigate the existence/ non-existence of solution of the following problem:
\begin{equation}
\opmb{u}{\O}+f(x,u)=0 \quad \text{ in } \quad \O
\end{equation}
where $f$ is of KPP type. We  show that the existence of a non trivial solution of \eqref{pev.eq.semilin} is governed by the sign of the principal eigenvalue of the following operator $ \mb{\O}+f_{u}(x,0).$ Moreover, when a non  trivial solution exists, then it is unique.

To show the existence/ non existence of solutions of \eqref{pev.eq.semilin} and their properties,   we follow and adapt  the arguments developed  in \cite{BHR1,BHRo,CDM2}.

\subsection{Existence of a non trivial solution}~
 
 Let us assume that  
$$
\lambda_p(\mb{\O}+f_u(x,0))<0.
$$
Then we we will show that there exists a non trivial solution to \eqref{pev.eq.semilin}.

Before going to the construction of a non trivial solution,  let us first define some quantities. First let us denote $a(x):=f_u(x,0)-b(x)$ and $\sigma:=\sup_\O a(x)$. 
 Observe that with this notation, we have $\lambda_p(\mb{\O}+f_u(x,0))=\lambda_p(\lb{\O}+a(x))$.

 From the definition of $\sigma$ there exists a sequence of points  $(x_n)_{n\in \N}$ such that $x_n \in \O$ and $ |\sigma -a(x_n)|\le \frac{1}{n}$.

Then by continuity of  $a(x)$, for each $n$ there exists $\eta_n$ such that for all $x \in B_{\eta_n}(x_n)$ we have 
$ |\sigma -a(x)|\le \frac{2}{n}$.

Now let us consider a sequence of real numbers $(\eps_n)_{n\in\N}$ which converges to zero such that  $\eps_n\le\frac{ \eta_n}{2}$.
 
Next, let $(\chi_n)_{n\in \N}$ be  the following sequence of  cut-off" functions :
 $\chi_n(x):=\chi(\frac{\|x-x_n\|}{\eps_n})$ where $\chi$ is a smooth function such that $0\le \chi \le 1$, $\chi(x) = 0$ for
$|x| \ge 2$ and $\chi(x) = 1$ for $|x| \le 1$.

Finally, let us consider the following sequence of continuous functions $(a_n)_{n\in \N}$, defined by   $a_n(x):=sup\{a(x),\sigma\chi_n\}$.
Observe that by construction the sequence $(a_n)_{n\in \N}$  is such that   $\|a(x)-a_n(x)\|_{\infty}\to 0$.

Let us now proceed to the construction of a non trivial solution.

By construction, for each $n$, the  function $a_n$ satisfies  $\sup_{\O} a_n=\sigma$ and  $a_n\equiv \sigma  $ in $B_{\frac{\eps_n}{2}}(x_n)$.  Therefore, the sequence $a_n$ satisfies $\frac{1}{\sigma-a_n}\not \in L^{1}_{loc}(\O)$ and by Theorem \ref{pev.th1} there exists a principal eigenpair $(\lambda_p^n,\phi_n)$ solution of the eigenvalue problem:
$$\oplb{\phi}{\O}+a_n(x)\phi+\lambda \phi=0,$$
such that $\phi_n\in L^{\infty}(\O)\cap C(\O)$.

Next,  using that $\|a_n(x) -a(x)\|_{\infty} \to 0$ as $n \to \infty$, from (iii) of the Proposition \ref{pev.prop1} it follows that for $n$ big enough,  say $n\ge n_0$, we have
 $$
\lambda_p^n< \frac{\lambda_p(\lb{\O}+a(x))}{2}<0.
$$
Moreover, by choosing $n_0$ bigger if necessary, we achieve for $n\ge n_0$
$$\lambda_p^n +\|a_n(x) -a(x)\|_{\infty} \le   \frac{\lambda_p(\lb{\O}+a(x))}{4}.$$

Let us now compute $ \opmb{\eps\phi_n}{\O} +f(x,\eps\phi_n)$.  
For $n\ge n_0$, we have
\begin{align*}
\opmb{\eps\phi_n}{\O} +f(x,\eps\phi_n)&=f(x,\eps\phi_n)-(b(x)+a_n(x))\eps\phi_n-\eps\lambda_p^n  \phi_n\\
&= \big(f_u(x,0)- (a_n(x)+b(x))\big)\eps\phi_n-\eps\lambda_p^n  \phi_n+o(\eps\phi_n)\\
&\ge \big(-\| a(x)- a_n(x)\|_{\infty} -\lambda_p^n\big)\eps\phi_n+o(\eps\phi_n)\\
& \ge -\frac{\lambda_p(\mb{\O}+f_u(x,0))}{4}\eps\phi_n+o(\eps\phi_n)>0.
\end{align*}
Therefore, for $\eps>0$ sufficiently small and $n$ big enough, $\eps\phi_n$ is a 
subsolution of \eqref{pev.eq.semilin}. By definition of $f$, any
large enough constant $M$  is a  supersolution of \eqref{pev.eq.semilin}.  By choosing $M$ so large that $\eps\phi_n\le M$ and using
a basic iterative scheme we obtain  the existence of a positive non trivial
 solution $u$ of \eqref{pev.eq.semilin}.
 

\subsection{Non-existence of positive bounded solutions}~

Let now turn our attention to the non-existence result. Let us prove that  when  $\lambda_p(\mb{\O}+f_u(x,0))\ge 0$ then there exists no non trivial solution to \eqref{pev.eq.semilin}.

Assume by contradiction that $\lambda_p(\mb{\O}+f_u(x,0))\ge 0 $ and  there  exists a positive bounded solution  $u$ to equation \eqref{pev.eq.semilin}.
  
  Obviously,  since $u$ is  nonnegative  and bounded, using \eqref{pev.eq.semilin}  we have for all $x\in \bar \O$
\begin{equation}
 0\le \oplb{u}{\O}=(b(x)-\frac{f(x,u)}{u})u.
\end{equation} 

Let us denote  $h(x):=\oplb{u}{\O}$. By construction, $h$ is a nonnegative continuous function in $\bar \O$. Therefore, since $\bar \O$ is compact,  $h$ achieves at some point  $x_0 \in \bar \O$ a nonnegative minimum.  
A short argument show that $h(x_0)>0$. Indeed, otherwise we have
$$ \int_{\O}J\left(\frac{x_0-y}{g(y)}\right)\frac{u(y)}{g^{n}(y)}\, dy =0.$$
Thus, since $J, g$ and $u$ are nonnegative quantities, from the above equality we deduce that  $u(y)= 0$ for almost every $y \in \{z\in \bar \O|\, \frac{x_0-z}{g(z)} \in supp(J)  \}$. By iterating this argument and using the assumption $J(0)>0$, we can show  that $u(y)=0$    for almost every $y\in \bar \O$, which implies that $u\equiv 0$ since $u$ is  continuous.
 
As a consequence  $\inf_{x\in\O}(b(x)-\frac{f(x,u)}{u})\ge \delta $ for some $\delta>0$  and  there exists a positive constant $c_0$ so that $u>c_0$ in $\bar \O$. 
 From the monotone properties of $f(x,.)$, we deduce that $\frac{f(x,u)}{u}\le \frac{f(x,c_0)}{c_0}<f_{u}(x,0)$.
 Let us now denote $\gamma(x)=\frac{f(x,c_0)}{c_0} -b(x)$. By construction, we have $\gamma(x)< a(x)$ and therefore by (ii) of Proposition \ref{pev.prop1}, $$\lambda_p(\lb{\O}+\gamma(x))>\lambda_{p}(\lb{\O}+a(x))\ge 0.$$    
Moreover, since $u$ is a solution of \eqref{pev.eq.semilin}, we have
$$\oplb{u}{\O}+\gamma(x)u\ge \opmb{u}{\O}+f(x,u)=0.$$

By definition of $\lambda_p(\lb{\O}+\gamma(x))$,  for all positive $\lambda <\lambda_p(\lb{\O}+\gamma(x))$ there exists a positive continuous function $\phi_\lambda$ such that 
$$\oplb{\phi_\lambda}{\O}+\gamma(x)\phi_\lambda\le -\lambda \phi_\lambda \le 0.$$ 

Arguing as  above, we can see that  $\phi_\lambda\ge\delta$ for some positive $\delta$.  Let us 
define the following quantity
$$\tau^*:=\inf\{\tau>0|u\le\tau\phi_\lambda\}.$$

Obviously, we end the proof of the theorem by proving  that $\tau^*=0$.
Assume that $\tau^*>0$. Then by
definition of $\tau^*$,  there exists $x_0\in \bar \O$ such that
$\tau^*\phi_p(x_0)=u(x_0)>0$.
At this point $x_0$, we have,
$$
0\le
\oplb{w}{\O}(x_0)=\oplb{(\tau^*\phi_\lambda-u)}{\O}(x_0)\le
0.
$$
Therefore, since $w\ge 0$, using a similar argumentation as above,  we have  $w(y)=0$ for almost every $y \in\bar \O$.
Thus, we end up with   $\tau^*\phi_1\equiv u$ and  we get the following  contradiction,
$$
0\le  \oplb{ u}{\O}+\gamma(x)u=\oplb{ \tau^* \phi_\lambda}{\O}+\gamma(x)\tau^*\phi_\lambda <  0.
$$
Hence $\tau^*=0$.

\fdem

\subsection{Uniqueness of the solution}~

Lastly,  we show that when a solution of  \eqref{pev.eq.semilin} exists then it is unique.
The proof of the uniqueness of the solution  is obtained as follows. 

Let $u$ and $v$ be two nonnegative bounded solution of \eqref{pev.eq.semilin}. Arguing as in the  above sub-section, we see that  there exists two positive constants $c_0$ and $c_1$ such that 
\begin{align*}
&u\ge c_0 \quad \text{ in } \bar \O\\
 &v\ge c_1 \quad \text{ in } \bar \O.
 \end{align*}
 
Since $u$ and $v$ are bounded and  strictly positive, the following quantity is well defined
$$ \gamma^*:=\inf\{\gamma>0\, |\, \gamma u\ge v \}.$$
We claim that $\gamma^*\le 1$. Indeed, assume by contradiction that $\gamma^*>1$. 
From \eqref{pev.eq.semilin} we see that  
\begin{align}
\opmb{\gamma^*u}{\O}+f(x,\gamma^*u)&= f(x,\gamma^* u)-\gamma^*f(x,u)\\
 &=\gamma^* u\left(\frac{f(x,\gamma^* u)}{\gamma^* u}-\frac{f(x,u)}{u}\right)\le 0\label{pev.eq.uni1}
 \end{align}

Now, by definition of $\gamma^*$, there exists $x_0\in \bar\O$
so that $\gamma u(x_0)= v(x_0)$ and 
from \eqref{pev.eq.semilin} we can easily see that 
\begin{equation}
\opmb{\gamma^*u}{\O}(x_0)+f(x,\gamma^*u(x_0))= \oplb{\gamma^* u-v}{\O}\ge 0.\label{pev.eq.uni2}
\end{equation}
From \eqref{pev.eq.uni1} and \eqref{pev.eq.uni2} we deduce that  
$$\oplb{\gamma^* u-v}{\O}(x_0)= 0.$$
Therefore,  arguing as in the above sub-section it follows  that $\gamma^* u =v$. Using now \eqref{pev.eq.uni1}, we deduce that 
 
$$0=\opmb{v}{\O}+f(x,v)=\opmb{\gamma^*u}{\O}+f(x,\gamma^*u) =\gamma^* u\left(\frac{f(x,\gamma^* u)}{\gamma^* u}-\frac{f(x,u)}{u}\right) \le 0,$$
 which implies that for all $x\in \O$ $f(x,\gamma^*u)\equiv f(x,u)$. This later is impossible since $\gamma^*>1$.
Hence, $\gamma^*\le 1$ and as a consequence $u\ge v$.

Observe that the role of $u$ and $v$ can be interchanged in the above argumentation. So we also have $v\ge u$, which shows the uniqueness of the solution.

\fdem
\section{Asymptotic Behaviour of the solution of  \eqref{pev.eq.semilin-para}}
Lastly, in this section, we prove the Theorem \ref{pev.thab} which establishes the asymptotic behaviour of the solution of 
	\begin{align*}
&\frac{\partial u}{\partial t}=\opmb{u}{\O} +f(x,u) \quad \text{ in }\quad \R^+\times\O.\\
&u(0,x)=u_0(x)\quad \text{ in }\quad \O
\end{align*}
\dem{ Proof of Theorem \ref{pev.thab}:}

The existence of a solution defined for all time $t$ follows from a standard argument and will not be exposed. Moreover, since $u_0\ge 0$ and $u_0\not\equiv 0$, using the parabolic maximum principle, there exists a positive constant  $\delta$ such that  $u(1,x)>\delta$ in $\bar \O$.
Let us first assume that $\lambda_p<0$. By following the argument developed  in above section, we can construct a bounded continuous function $\psi$ so  that  $\eps \psi$ is a subsolution of \eqref{pev.eq.semilin-para} for $\eps$ small enough. Since, $u(1,x)\ge \delta$ and $\psi$ is bounded, by choosing $\eps$ smaller if necessary we achieves also that $\eps \psi \le u(1,x)$.
Now, let us denote $\under \Psi(x,t)$  the solution of evolution problem \eqref{pev.eq.semilin-para} with initial datum $\eps\psi$. By construction, using a standard argument,  $\under\Psi(t,x)$ is a non-decreasing function of the time and $\Psi(t,x)\le u(t+1,x)$.
On the other hand, since for $M$ big enough $M$ is a supersolution of  \eqref{pev.eq.semilin-para} and $u_0$ is bounded, we have also  
$u(t,x)\le\bar \Psi(t,x), $ where   $\bar \Psi(x,t)$ denotes the solution of evolution problem \eqref{pev.eq.semilin-para} with initial datum $\bar \Psi(0,x)=M\ge u_0$. A standard argument using the parabolic comparison principle shows that $\bar \Psi$ is a non-increasing function of $t$. 
Thus we have for all time $t$ 
$$\eps\psi\le \under\Psi(t,x)\le u(t+1,x)\le \bar \Psi(t+1,x).$$
Since $\under\Psi(t,x)$ (respectively $ \bar \Psi(t,x)$)  is an  uniformly bounded monotonic function of $t$, $\under \Psi$ (resp. $\bar \Psi$) converges pointwise  to $\under p$ (resp. $\bar p$) which is a solution of  \eqref{pev.eq.semilin}. From  $\under\Psi(t,x) \not\equiv 0$, using the uniqueness of a non-trivial solution (Theorem \ref{pev.th3}), we deduce that   $\under p\equiv \bar p \not \equiv 0$ and therefore,
$u(x,t)\to p$ pointwise in $\O$, where $p$ denotes the unique non trivial solution of \eqref{pev.eq.semilin}.

In the other case, when $\lambda_p\ge 0$ we argue as follows.
As above, we have $0\le u(t,x)\le\bar \Psi(t,x) $ and $\bar \Psi $ converges pointwise to   $\bar p$ a  solution of  \eqref{pev.eq.semilin}. 
By Theorem \ref{pev.th3}in this situation we have $\bar p\equiv 0$, hence $u(x,t)\to 0$ pointwise in  $\O$. 

\fdem

\begin{remark}
Note that the above analyse will hold for more general kernel non negative kernel $k(x,y)$ that satisfies $\tilde H2$, i.e.
$$
\exists\, c_0>0,\ \eps_0 >0\; \text{ such that } \min_{x\in\O}\left(\min_{y\in B(x,\eps_0)}k(x,y)\right)>c_0.
$$ 
\end{remark}
\appendix
\section{}
In this appendix, we first prove the Proposition \ref{pev.prop1}. Then we recall the method of sub and supersolution to obtain solution of the semilinear problem :
\begin{equation}\label{pev-eq-semilinapp}
\opmb{u}{\O}=f(x,u) \quad \text{ in } \quad \O.
\end{equation}

Before going to the proof of the Proposition \ref{pev.prop1}, let us show that $\lambda_p(\lb{\O}+a(x))$ is well defined.
Let us first show that the set $\Lambda:=\{\lambda \,|\, \exists \, \phi \in C(\O), \phi>0 \; \text{ such that } \; \oplb{\phi}{\O}+\lambda \phi \le 0\}$ is non-empty. 
Indeed, as observed in \cite{Co5} (Theorem 1.8), for $\O,J,g$ and $a$ satisfying the assumptions (H1-H4)  there exists a continuous positive function $\psi$ satisfying 
$$ \int_{\O}J\left(\frac{x-y}{g(y)}\right)\frac{\psi(y)}{g^n(y)}\, dy= c(x)\psi(x),$$
 where $c(x)$ is defined by 
 \begin{equation*}
c(x):=\left\{
\begin{array}{l}
1 \quad  \text{ if } \quad x \in \{x\in \bar\O\,|\, g(x)=0\}\\
\int_{\O}J\left(\frac{y-x}{g(x)}\right)\frac{dy}{g^n(x)}\quad \text{ otherwise. }
\end{array}\right.
\end{equation*}
Obviously $c(x)\in L^{\infty}$ and  for any $ \lambda\le\-(| a\|_{\infty}+\| c\|_{\infty}) $ we have 
\begin{align*}
\oplb{\psi}{\O}+(a(x)+\lambda)\psi&=(a(x)+c(x)+\lambda)\psi\\
&\le( a(x)+c(x)-\| a\|_{\infty}-\| c\|_{\infty})\psi \le 0.
\end{align*}
Therefore, the set $\Lambda$ is non-empty. 

Observe now that since $J,g$ are nonnegative functions and  $a(x)\in L^{\infty}$, for any continuous positive function $\phi$ we have
$$ \oplb{\phi}{\O}+(a(x)+\|a(x)\|_{\infty})\phi \ge 0.$$
Therefore, the set $\Lambda$ has an upper bound and $\lambda_p$ is well defined.

Let us now  prove  the Proposition \ref{pev.prop1}.

\dem{Proof of the Proposition \ref{pev.prop1} :}
(i) easily follows from the definition of  $\lambda_p$. First, let us observe that to obtain 
$$   \lambda_p(\lb{\O_2}+a(x))\le \lambda_p(\lb{\O_1}+a(x))$$
it is sufficient to prove the inequality 
$$ \lambda\le \lambda_p(\lb{\O_1}+a(x))$$ for any $\lambda <\lambda_p(\lb{\O_2}+a(x))$.

 Let us fix $\lambda< \lambda_p(\lb{\O_2}+a(x))$. Then by definition of $\lambda_p(\lb{\O_2}+a(x))$ 
 there exists  a positive function $ \phi \in C(\O_2)$ such that
 $$\oplb{\phi}{\O_2}+(a(x)+\lambda)\phi\le 0. $$
Since $\O_1 \subset \O_2$, an easy computation shows that 
$$
\oplb{\phi}{\O_1}+(a(x)+\lambda)\phi \le \oplb{\phi}{\O_2}+(a(x)+\lambda)\phi \le 0
$$
Therefore, by definition of $\lambda_p(\lb{\O_1}+a(x))$ we have $\lambda \le \lambda_p(\lb{\O_1}+a(x))$. Hence,
$\lambda_p(\lb{\O_2}+a(x)) \le \lambda_p(\lb{\O_1}+a(x))$.
\smallskip

To show (ii), we argue as above. By definition of $\lambda_p(\lb{\O}+a_1(x)) $ for any $\lambda <\lambda_p(\lb{\O}+a_1(x)) $   there exists a positive  $\phi\in C(\O)$ such that
$$ \oplb{\phi}{\O}+(a_1(x)+\lambda)\phi\le 0$$ 
and we  have 
$$ \oplb{\phi}{\O}+(a_2(x)+\lambda)\phi \le  \oplb{\phi}{\O}+(a_1(x)+\lambda)\phi\le 0.$$
Therefore $\lambda\le \lambda_p(\lb{\O}+a_2(x)) $. Hence (ii) holds true.
\smallskip

Let us now prove (iii). Again we fix $\lambda<\lambda_p(\lb{\O}+a(x))$. For this $\lambda$, there exists $\phi \in C(\O)$, $\phi>0$ such that
\begin{equation}
\oplb{\phi}{\O}+(a(x)+\lambda)\phi \le 0.\label{pev.eq.prop1}
\end{equation}

An easy computation shows that we rewrite the above equation the following way:
\begin{align*}
\oplb{\phi}{\O}+(a(x)+\lambda)\phi &= \oplb{\phi}{\O}+(b(x)+\lambda)\phi +(a(x)-b(x))\phi\\
&\ge \oplb{\phi}{\O}+(b(x)+\lambda-\| a(x)-b(x)\|_{\infty})\phi
\end{align*}
Using that $(\lambda,\phi)$ satisfies \eqref{pev.eq.prop1}, it follows that 
 $$\oplb{\phi}{\O}+(b(x)+\lambda-\| a(x)-b(x)\|_{\infty})\phi\le 0.$$
Therefore, $\lambda-\| a(x)-b(x)\|_{\infty}\le \lambda_p(\lb{\O}+b(x))$ and we have
$$\lambda\le \lambda_p(\lb{\O}+b(x))+\| a(x)-b(x)\|_{\infty}.$$

The above computation being valid for any $\lambda<\lambda_p(\lb{\O}+a(x))$, we end up with
$$ \lambda_p(\lb{\O}+a(x))-\lambda_p(\lb{\O}+b(x))\le \| a(x)-b(x)\|_{\infty}. $$

Note that  the role of $a(x)$ and $b(x)$ can be interchanged in the above argumentation. So, we also have
  $$ \lambda_p(\lb{\O}+b(x))-\lambda_p(\lb{\O}+a(x))\le \| a(x)-b(x)\|_{\infty}. $$
  Hence
$$| \lambda_p(\lb{\O}+a(x))-\lambda_p(\lb{\O}+b(x))|\le \| a(x)-b(x)\|_{\infty},$$
which proves (iii).

The proof of (iv) being   similar to the proof of (ii), it will be omitted. 
\fdem
  
Before recalling the sub/supersolution method, let us introduce some definitions and notations.
We call a  bounded continuous function $\bar u$ (resp. $ \under{u}$) a super-solution (resp. a sub-solution) if  $\bar u$ (resp. $ \under{u}$)   satisfies the following inequalities: 
\begin{equation}
\opmb{u}{\O}\le(\ge) f(x,u) \quad \text{ in } \quad \O.
\end{equation}
Let us now state the  Theorem.
\begin{theorem}
Assume $f(x,.)$ is a Lipschitz function uniformly in $x$ and let $\bar u$ and $\under{u}$ be respectively a supersolution  and a subsolution of \eqref{pev-eq-semilinapp} continuous up to the boundary. Assume further that $\under{u}\le \bar u$. Then there exists a solution $u\in C(\bar\O)$ solution of \eqref{pev-eq-semilinapp} satisfying $\under{u}\le u\le \bar u$     
\end{theorem}

\dem{Proof:}

Let us  first choose  $k>|\lambda_p(\mb{\O})|$ big enough such that the function $-ks+f(x,s)$ is a decreasing function of $s$ uniformly in $x$. 
We can increase further $k$ if necessary to ensure that $k\in \rho(\mb{\O}) $, where $\rho(\mb{\O}) $ denotes the resolvent of the operator $\mb{\O}$.

Note that by this choice of $k$, by Theorem \ref{pev.th2} the operator $\mb{\O}- k$ satisfies a comparison principle. 

Now, let $u_1$ be the solution of the following linear problem 
\begin{equation}\label{pev-eq-induc1}
\opmb{u_1}{\O} -ku_1=-k\under{u}+f(x,\under{u}) \quad\text{ in }\quad \O.
\end{equation}
$u_1$ always exists, since by construction the continuous operator $\mb{\O}-k$ is invertible.
We claim that $\under{u}\le u_1 \le \bar u$. Indeed, since $\under{u}$ and $\bar u$ are respectively a sub- and super-solution of \eqref{pev-eq-semilinapp}, we have 
\begin{align*}
&\opmb{u_1-\under{u}}{\O}-k(u_1-\under{u})\le 0 \quad\text{ in }\quad \O\\
&\opmb{u_1-\bar u}{\O}-k(u_1-\bar u)\ge -k(\under{u}-\bar u) +f(x,\under{u})-f(x,\bar u)\ge 0 \quad\text{ in }\quad \O.
\end{align*}
So, the inequality $\under{u}\le u_1 \le \bar u$ follows from the comparison principle satisfied by the operator $\mb{\O}-k$.
Now let  $u_2$ be  the solution of \eqref{pev-eq-induc1} with $u_1$ instead of $\under{u}$.  From the monotonicity of $-ks+f(x,s)$ and  using the  comparison principle, we have $\under{u}\le u_1\le u_2 \le \bar u$. By induction, we can construct an increasing  sequence of function $(u_n)_{n\in\N}$ satisfying  $\under{u}\le u_n \le \bar u$ and 
\begin{equation}\label{pev-eq-induc}
\opmb{u_{n+1}}{\O} -ku_{n+1}=-ku_n+f(x,u_n) \quad\text{ in }\quad \O.
\end{equation}
Since the sequence is increasing and bounded, $u^-(x):=sup_{n\in \N}u_n(x)$ is well defined.  Moreover, passing to the limit in the equation \eqref{pev-eq-induc} using  Lebesgue's Theorem it follows that  $u^-$ is a solution of \eqref{pev-eq-semilinapp}.  
\fdem

\section*{}
\bibliographystyle{plain}
\bibliography{pev.bib}

\def\cprime{$'$}
\begin{thebibliography}{10}

\bibitem{AB}
G.~Alberti and G.~Bellettini.
\newblock A nonlocal anisotropic model for phase transitions. {I}. {T}he
  optimal profile problem.
\newblock {\em Math. Ann.}, 310(3):527--560, 1998.

\bibitem{BC1}
P.~W. Bates and A.~Chmaj.
\newblock An integrodifferential model for phase transitions: stationary
  solutions in higher space dimensions.
\newblock {\em J. Statist. Phys.}, 95(5-6):1119--1139, 1999.

\bibitem{BFRW}
P.~W. Bates, P.~C. Fife, X.~Ren, and X.~Wang.
\newblock Traveling waves in a convolution model for phase transitions.
\newblock {\em Arch. Rational Mech. Anal.}, 138(2):105--136, 1997.

\bibitem{BZ}
P.~W. Bates and G.~Zhao.
\newblock Existence, uniqueness and stability of the stationary solution to a
  nonlocal evolution equation arising in population dispersal.
\newblock {\em J. Math. Anal. Appl.}, 332(1):428--440, 2007.

\bibitem{BHR1}
H.~Berestycki, F.~Hamel, and L.~Roques.
\newblock Analysis of the periodically fragmented environment model. {I}.
  {S}pecies persistence.
\newblock {\em J. Math. Biol.}, 51(1):75--113, 2005.

\bibitem{BHRo}
H.~Berestycki, F.~Hamel, and L.~Rossi.
\newblock Liouville-type results for semilinear elliptic equations in unbounded
  domains.
\newblock {\em Ann. Mat. Pura Appl. (4)}, 186(3):469--507, 2007.

\bibitem{BNV}
H.~Berestycki, L.~Nirenberg, and S.~R.~S. Varadhan.
\newblock The principal eigenvalue and maximum principle for second-order
  elliptic operators in general domains.
\newblock {\em Comm. Pure Appl. Math.}, 47(1):47--92, 1994.

\bibitem{browder}
F.~E. Browder.
\newblock On the spectral theory of elliptic differential operators. {I}.
\newblock {\em Math. Ann.}, 142:22--130, 1960/1961.

\bibitem{CMS}
M.~L. Cain, B.~G. Milligan, and A.~E. Strand.
\newblock {Long-distance seed dispersal in plant populations}.
\newblock {\em Am. J. Bot.}, 87(9):1217--1227, 2000.

\bibitem{CCR}
E.~Chasseigne, M.~Chaves, and J.~D. Rossi.
\newblock Asymptotic behavior for nonlocal diffusion equations.
\newblock {\em J. Math. Pures Appl. (9)}, 86(3):271--291, 2006.

\bibitem{Ch}
X.~Chen.
\newblock Existence, uniqueness, and asymptotic stability of traveling waves in
  nonlocal evolution equations.
\newblock {\em Adv. Differential Equations}, 2(1):125--160, 1997.

\bibitem{CR}
A.~Chmaj and X.~Ren.
\newblock The nonlocal bistable equation: stationary solutions on a bounded
  interval.
\newblock {\em Electron. J. Differential Equations}, pages No. 02, 12 pp.
  (electronic), 2002.

\bibitem{Cl}
J.~S. Clark.
\newblock Why trees migrate so fast: Confronting theory with dispersal biology
  and the paleorecord.
\newblock {\em The American Naturalist}, 152(2):204--224, 1998.

\bibitem{CCEM}
C.~Cort{\'a}zar, J.~Coville, M.~Elgueta, and S.~Mart{\'{\i}}nez.
\newblock A nonlocal inhomogeneous dispersal process.
\newblock {\em J. Differential Equations}, 241(2):332--358, 2007.

\bibitem{CER}
C.~Cort{\'a}zar, M.~Elgueta, and J.~D. Rossi.
\newblock A nonlocal diffusion equation whose solutions develop a free
  boundary.
\newblock {\em Ann. Henri Poincar\'e}, 6(2):269--281, 2005.

\bibitem{Co2}
J.~Coville.
\newblock On uniqueness and monotonicity of solutions of non-local reaction
  diffusion equation.
\newblock {\em Ann. Mat. Pura Appl. (4)}, 185(3):461--485, 2006.

\bibitem{AC}
J.~Coville.
\newblock Remarks on the strong maximum principle for nonlocal operators.
\newblock {\em Electron. J. Differential Equations}, pages No. 66, 10, 2008.

\bibitem{Co5}
J.~Coville.
\newblock Harnack's inequality for some nonlocal equations and application.
\newblock {\em Preprint du MPI}, Feb. 2008.

\bibitem{Co4}
J.~Coville.
\newblock Travelling fronts in asymmetric nonlocal reaction diffusion equation:
  The bistable and ignition case.
\newblock {\em Preprint du CMM}, July. 2006.

\bibitem{CDM2}
J.~Coville, J.~D\'{a}vila, and S.~Mart\'{\i}nez.
\newblock Existence and uniqueness of solutions to a nonlocal equation with
  monostable nonlinearity.
\newblock {\em SIAM Journal on Mathematical Analysis}, 39(5):1693--1709, 2008.

\bibitem{CDM1}
J.~Coville, J.~D{\'a}vila, and S.~Mart{\'{\i}}nez.
\newblock Nonlocal anisotropic dispersal with monostable nonlinearity.
\newblock {\em J. Differential Equations}, 244(12):3080--3118, 2008.

\bibitem{CD1}
J.~Coville and L.~Dupaigne.
\newblock Propagation speed of travelling fronts in non local
  reaction-diffusion equations.
\newblock {\em Nonlinear Anal.}, 60(5):797--819, 2005.

\bibitem{CD2}
J.~Coville and L.~Dupaigne.
\newblock On a non-local equation arising in population dynamics.
\newblock {\em Proc. Roy. Soc. Edinburgh Sect. A}, 137(4):727--755, 2007.

\bibitem{DGP}
A.~De~Masi, T.~Gobron, and E.~Presutti.
\newblock Travelling fronts in non-local evolution equations.
\newblock {\em Arch. Rational Mech. Anal.}, 132(2):143--205, 1995.

\bibitem{DK}
C.~Deveaux and E.~Klein.
\newblock Estimation de la dispersion de pollen \`a longue distance \`a
  l'echelle d'un paysage agicole : une approche exp\'erimentale.
\newblock {\em Publication du Laboratoire Ecologie, Syst\`ematique et
  Evolution}, 2004.

\bibitem{DV}
M.~D. Donsker and S.~R.~S. Varadhan.
\newblock On a variational formula for the principal eigenvalue for operators
  with maximum principle.
\newblock {\em Proc. Nat. Acad. Sci. U.S.A.}, 72:780--783, 1975.

\bibitem{EPS}
D.~E. Edmunds, A.~J.~B. Potter, and C.~A. Stuart.
\newblock Non-compact positive operators.
\newblock {\em Proc. Roy. Soc. London Ser. A}, 328(1572):67--81, 1972.

\bibitem{E}
L.~C. Evans.
\newblock {\em Partial differential equations}, volume~19 of {\em Graduate
  Studies in Mathematics}.
\newblock American Mathematical Society, Providence, RI, 1998.

\bibitem{F1}
P.~C. Fife.
\newblock {\em Mathematical aspects of reacting and diffusing systems},
  volume~28 of {\em Lecture Notes in Biomathematics}.
\newblock Springer-Verlag, Berlin, 1979.

\bibitem{F2}
P.~C. Fife.
\newblock An integrodifferential analog of semilinear parabolic {PDE}s.
\newblock In {\em Partial differential equations and applications}, volume 177
  of {\em Lecture Notes in Pure and Appl. Math.}, pages 137--145. Dekker, New
  York, 1996.

\bibitem{GR}
J.~Garc\'ia-Meli\'an and J.~D. Rossi.
\newblock On the principal eigenvalue of some nonlocal diffusion problems.
\newblock {\em Journal of Differential Equations}, 246(1):21 -- 38, 2009.

\bibitem{HMMV}
V.~Hutson, S.~Martinez, K.~Mischaikow, and G.~T. Vickers.
\newblock The evolution of dispersal.
\newblock {\em J. Math. Biol.}, 47(6):483--517, 2003.

\bibitem{KM}
J.~Medlock and M.~Kot.
\newblock Spreading disease: integro-differential equations old and new.
\newblock {\em Math. Biosci.}, 184(2):201--222, 2003.

\bibitem{M}
J.~D. Murray.
\newblock {\em Mathematical biology}, volume~19 of {\em Biomathematics}.
\newblock Springer-Verlag, Berlin, second edition, 1993.

\bibitem{nussbaum}
R.~D. Nussbaum.
\newblock The radius of the essential spectrum.
\newblock {\em Duke Math. J.}, 37:473--478, 1970.

\bibitem{NP}
R.~D. Nussbaum and Y.~Pinchover.
\newblock On variational principles for the generalized principal eigenvalue of
  second order elliptic operators and some applications.
\newblock {\em J. Anal. Math.}, 59:161--177, 1992.
\newblock Festschrift on the occasion of the 70th birthday of Shmuel Agmon.

\bibitem{PW}
M.~H. Protter and H.~F. Weinberger.
\newblock {\em Maximum principles in differential equations}.
\newblock Prentice-Hall Inc., Englewood Cliffs, N.J., 1967.

\bibitem{P}
C.~Pucci.
\newblock Maximum and minimum first eigenvalues for a class of elliptic
  operators.
\newblock {\em Proc. Amer. Math. Soc.}, 17:788--795, 1966.

\bibitem{SSN}
F.~M. Schurr, O.~Steinitz, and R.~Nathan.
\newblock Plant fecundity and seed dispersal in spatially heterogeneous
  environments: models, mechanisms and estimation.
\newblock {\em J. Ecol.}, 96(4):628--641, 2008.

\bibitem{zeidler}
E.~Zeidler.
\newblock {\em Nonlinear functional analysis and its applications. {I}}.
\newblock Springer-Verlag, New York, 1986.
\newblock Fixed-point theorems, Translated from the German by Peter R. Wadsack.

\end{thebibliography}

\end{document}